**Torus bundles and the group cohomology of $GL(N, \mathbb{Z})$**

by

**Jean-Michel Bismut**

Départment de Mathématique

Université Paris-Sud

91405 - Orsay

France

and

**John Lott**

Department of Mathematics

University of Michigan

Ann Arbor, MI 48109-1003

USA

lott@math.lsa.umich.edu

**Abstract**

We prove the vanishing of a certain characteristic class of flat vector bundles when the structure groups of the bundles are contained in $GL(N, \mathbb{Z})$. We do so by explicitly writing the characteristic class as an exact form on the base of the bundle.

CONTENTS

In this paper we consider certain characteristic classes of flat complex vector bundles, which are known in algebraic K-theory as the Borel regulator classes. We prove that if the structure group of a rank-$N$ vector bundle is contained in $GL(N, \mathbb{Z})$ with $N$ odd then the Borel class of degree $2N - 1$ vanishes. Our proof is analytic in nature and is a special, but interesting, case of a more general theorem concerning the direct images of flat vector bundles under smooth submersions [BLo].

The background to our result is the following. First, let $N$ be a positive even integer and let $E$ be a real oriented rank-$N$ vector bundle over a connected manifold $B$. Then the rational Euler class $\chi_{\mathbb{Q}}(E)$ is an element of $\mathrm{H}^N(B; \mathbb{Q})$. Sullivan showed that if $E$ is a flat vector bundle whose structure group is contained in $SL(N, \mathbb{Z})$, then $\chi_{\mathbb{Q}}(E) = 0$ [Su]. Let $\Lambda$ be the integer lattice in $E$. Then $M = E/\Lambda$ is the total space of a torus bundle over $B$. Sullivan's proof was by a simple topological argument involving this torus bundle.

Bismut and Cheeger observed that Sullivan's result follows from the Atiyah-Singer families index theorem, applied to the vertical signature operators on the torus bundle [BC2]. They also showed that one can write a certain differential-form representative of $\chi_{\mathbb{Q}}(E)$ explicitly as an exact form, by means of the so-called eta-form of [BC1].

In [BLo] we proved a real analog of the Riemann-Roch-Grothendieck (RRG) theorem. The geometric setup of our theorem involved a smooth fiber bundle $\pi : M \to B$ with closed connected fibers $Z$. If $F$ is a flat complex vector bundle on $M$ then it has a direct image $H\left(Z; F\big|_Z\right)$ which is an alternating sum of flat complex vector bundles on $B$, given by the cohomologies of the fibers $Z$ (with value in $F\big|_Z$). We defined certain characteristic classes of flat complex vector bundles and showed a relationship between the characteristic classes of $F$ and $H\left(Z; F\big|_Z\right)$. We actually proved a more refined statement at the level of differential forms on $B$, which involved a so-called analytic torsion form on $B$. These results are reviewed in Appendix A of the present paper.

The motivation of the present paper was to see what the RRG-type theorem of [BLo] says in the case of the torus bundle described above. We now state the precise results.

Let $B$ be a connected smooth manifold. Let $E$ be a complex rank-$N$ vector bundle over $B$ with a flat connection $\nabla^E$. We define certain characteristic classes $\{n_j^z(\nabla^E)\}_{j=1}^N$, with $n_j^z(\nabla^E) \in \mathrm{H}^{2j-1}(B; \mathbb{R})$. These classes are pulled back from universal classes $n_{j,\mathbb{C}^N}^z \in \mathrm{H}^{2j-1}(BGL(N, \mathbb{C})_\delta; \mathbb{R})$, where $\delta$ denotes the discrete topology on $GL(N, \mathbb{C})$. The classes $\{n^z_{j,\mathbb{C}^N}\}_{j=1}^N$ can be characterized by the fact that the continuous real-valued group cohomology of $GL(N,\mathbb{C})$ is an exterior algebra $\Lambda\left(n^z_{1,\mathbb{C}^N}, n^z_{2,\mathbb{C}^N}, \ldots, n^z_{N,\mathbb{C}^N}\right)$ [Bo]. They are stable classes in the sense that they come from classes in $\mathrm{H}^{2j-1}(BGL(\infty,\mathbb{C})_\delta; \mathbb{R})$, and give rise to the Borel regulators on the algebraic K-theory of number fields.

Our main result is the following :

**Theorem 0.1 :** *Let $N$ be a positive odd integer. Let $E$ be a flat complex rank-$N$ vector bundle over $B$ whose structure group is contained in $GL(N,\mathbb{Z})$. Then $n^z_N(\nabla^E)$ vanishes in $\mathrm{H}^{2N-1}(B;\mathbb{R})$.*

An equivalent formulation of the theorem is the following:

**Theorem 0.2 :** *Let $N$ be a positive odd integer. Let $i : GL(N,\mathbb{Z}) \to GL(N,\mathbb{C})$ be the natural inclusion. Then $i^*\left(n^z_{N,\mathbb{C}^N}\right)$ vanishes in the group cohomology $\mathrm{H}^{2N-1}(GL(N,\mathbb{Z}); \mathbb{R})$.*

Borel showed that for $j > 1$, the classes $i^*\left(n^z_{j,\mathbb{C}^N}\right)$ are nonzero if $N$ is sufficiently large compared to $j$ [Bo]. Ronnie Lee informs us that he showed many years ago in unpublished work that $i^*\left(n^z_{j,\mathbb{C}^N}\right)$ is nonzero if $1 < j < N$. Furthermore, Jens Franke informs us that Theorem 0.2 is a special case of general results of his on the group cohomology of arithmetic groups [F]. His method of proof seems to be completely different from ours.

Our proof of Theorem 0.1 parallels that of the corresponding proof of Sullivan's result in [BC2]. Given $t > 0$, in Section 2 we define a form $\delta_t$ on the total space of $E^*$ with the property that if $\lambda$ is a flat section of $E^*$ then $\lambda^* \delta_t$ is closed on $B$. Letting $\Lambda^*$ denote the dual lattice to $\Lambda$, we consider the closed form $\sum_{\mu \in \Lambda^*} \mu^* \delta_t$ on $B$. We show that its de Rham cohomology class $\left[\sum_{\mu \in \Lambda^*} \mu^* \delta_t\right]$ is independent of $t$. Taking $t \to 0$, one can see that $\left[\sum_{\mu \in \Lambda^*} \mu^* \delta_t\right]$ is a nonzero constant times $n^z_N(\nabla^E)$. We then define forms $\rho_t$ on the total space of $E$ which are again closed after being pulled back by flat sections, and which have the property that $\sum_{\mu \in \Lambda^*} \mu^* \delta_t$ is proportionate to $\sum_{m \in \Lambda} m^* \rho_t$. Taking $t \to \infty$, one sees that $\sum_{m \in \Lambda} m^* \rho_t \to 0$. This proves Theorem 0.1. By keeping track of the $t$-dependence, we find an explicit form on $B$, defined as a Dirichlet-type series, whose differential represents $n^z_N(\nabla^E)$.

The proof in Section 2 is completely elementary, in that it only uses Berezin integrals and the Poisson summation formula. However, it will undoubtedly seem unmotivated to the reader. The motivation comes from our previous work on direct images of flat vector bundles [BLo]. In Appendix A we give a summary of the results of [BLo] and show how they imply Theorems 0.1 and 0.2. This appendix can be read independently of Section 2. We then show in Appendix B that in the case of a flat torus bundle, one can use Fourier analysis to transcribe the results of Appendix A into the more elementary arguments of Section 2. This again parallels the arguments of [BC2].

An open question is whether there is a simple topological proof of Theorem 0.1 along the lines of Sullivan's proof of his result.

The paper is organized as follows. In Section 1 we define the relevant characteristic classes of flat vector bundles by means of differential-form representatives. Given a complex vector bundle $F$ with flat connection $\nabla^F$, a Hermitian metric $h^F$ on $F$ and a real invariant power series $P$ on the space of $(N \times N)$-complex matrices, we define a closed differential form $P^z(\nabla^F, h^F) \in \Omega^{\mathrm{odd}}(B)$. Its de Rham cohomology class $P^z(\nabla^F) \in \mathrm{H}^{\mathrm{odd}}(B; \mathbb{R})$ is independent of $h^F$. In the case of $P(A) = \mathrm{Tr}\left[A^j\right]$, we obtain the above-mentioned cohomology class $n_j^z(\nabla^F) \in \mathrm{H}^{2j-1}(B; \mathbb{R})$. We study the relationships between the classes and describe them in terms of the cohomology of the classifying space $BGL(N, \mathbb{C})_\delta$.

In Section 2 we start with a real oriented rank-$N$ vector bundle $\mathcal{V}$. We review Berezin integrals, the Thom form $\alpha_t \in \Omega^N(\mathcal{V})$ of Mathai-Quillen and its transgressing form $\beta_t \in \Omega^{N-1}(\mathcal{V})$. If $\mathcal{V}$ is flat, we define the form $\delta_t \in \Omega^{2N-1}(\mathcal{V})$ by a slight modification of the definition of $\alpha_t$. We also construct its transgressing form $\epsilon_t \in \Omega^{2N-2}(\mathcal{V})$. We then define an auxiliary form $\rho_t \in \Omega^{2N-1}(\mathcal{V})$ and its transgressing form $\sigma_t \in \Omega^{2N-2}(\mathcal{V})$. We prove Theorems 0.1 and 0.2 along the lines sketched above.

In Appendix A we give the relevant results from [BLo], without proof. We then look at the case of discrete groups acting on manifolds. Suppose that $\Gamma$ is a discrete group which has a CW classifying space $B\Gamma$ with a finite number of cells in each dimension. If $\Gamma$ acts smoothly on a compact manifold $Z$, we use the results of [BLo] to derive vanishing theorems in $\mathrm{H}^*(\Gamma; \mathbb{R})$. We then show that Theorem 0.2 arises as the special case of $GL(N, \mathbb{Z})$ acting on $T^N$.

In Appendix B we first apply the results of Appendix A to the total space of a flat vector bundle $\mathcal{V}$. Using the superconnection formalism of Appendix A, we define a form $\widetilde{\delta}_t \in \Omega^{2N-1}(\mathcal{V})$ and its transgressing form $\widetilde{\epsilon}_t \in \Omega^{2N-1}(\mathcal{V})$. We then show that $\widetilde{\delta}_t$ and $\widetilde{\epsilon}_t$ are essentially the same as the forms $\delta_t$ and $\epsilon_t$ of Section 2. Using Fourier analysis on the fibers of the torus bundle $M$, we show that the closed form $f\left(C'_t, h^W\right)$ of Appendix A reduces to $\sum_{\mu \in \Lambda^*} \mu^* \delta_t$, and that the transgressing form $\frac{1}{t} f^{\wedge} \left( C'_t, h^W \right)$ reduces to $\sum_{\mu \in \Lambda^*} \mu^* \epsilon_t$. This establishes the link between the results of Section 2 and Appendix A.

We thank Christophe Soulé and Toby Stafford for helpful discussions. The first author thanks the Institut Universitaire de France for its support. The second author thanks the NSF for its support.

## I - Characteristic classes of flat vector bundles

In this section we describe certain characteristic classes of flat vector bundles.

The section is organized as follows. In a) we establish our conventions. Given a rank-$N$ complex vector bundle $F$ on a base $B$, a flat connection $\nabla^F$ on $F$, a Hermitian metric $h^F$ on $F$ and an invariant power series $P$ on the space of complex $(N \times N)$-matrices, in b) we define a closed form $P^z(\nabla^F, h^F) \in \Omega^{\text{odd}}(B)$. Its de Rham cohomology class $P^z(\nabla^F)$ is independent of $h^F$. We describe relationships between the forms $P^z(\nabla^F, h^F)$ for different choices of $P$ and $F$. In c) we relate $P^z(\nabla^F)$ to the Borel regulator classes.

### a) Conventions

Except where otherwise indicated, we will take all vector spaces in this paper to be over $\mathbb{C}$. The covariant functor $\Lambda$ sends a vector space $V$ to its exterior algebra $\Lambda(V)$ and a linear map $T : V \to W$ to an algebra homomorphism $\Lambda(T) : \Lambda(V) \to \Lambda(W)$.

If $A$ is an $N \times N$ complex matrix, put

$$\begin{aligned}
(1.1) \quad c(A) &= \det(I + A), \\
\operatorname{ch}(A) &= \operatorname{Tr}\left[e^A\right], \\
\operatorname{Td}(A) &= \det\left(\frac{A}{I - e^{-A}}\right), \\
n_j(A) &= \operatorname{Tr}\left[A^j\right], j \in \mathbb{N}.
\end{aligned}$$

Let $\{c_j(A)\}_{j=1}^N$ be the symmetric functions of $A$, satisfying

$$(1.2) \qquad c(\lambda A) = 1 + \lambda c_1(A) + \ldots + \lambda^N c_N(A).$$

Let $B$ be a smooth connected manifold. If $E$ is a smooth vector bundle over $B$, we let $C^\infty(B; E)$ denote the smooth sections of $E$ and $L^2(B; E)$ denote the $L^2$ measurable sections of $E$. We let $\Lambda(T^*B)$ denote the complexified exterior bundle of $B$ and $\Omega(B)$ denote the space of smooth sections of $\Lambda(T^*B)$. We put $\Omega(B; E) = C^\infty(B; \Lambda(T^*B) \otimes E)$. We will say that a differential form is real if it can be written with real coefficients.

### b) Characteristic classes of flat vector bundles

Let $P(A)$ be an ad-invariant power series on the space of $N \times N$ complex matrices $A$. Then $P$ can be expressed as a power series in the variables $\{n_j(A)\}_{j=1}^N$. We will say that $P$ is real if $P$ has real coefficients in these variables.

Let $F$ be a complex rank-$N$ vector bundle on $B$, endowed with a flat connection $\nabla^F$. The antidual bundle $\overline{F}^*$ inherits a flat connection $\nabla^{\overline{F}^*}$. Let $h^F$ be a Hermitian metric on $F$. We do not require that $\nabla^F$ be compatible with $h^F$. The metric $h^F$ induces a Hermitian metric $h^{\overline{F}^*}$ on $\overline{F}^*$ and a $C^\infty(B)$-linear isometry

$$(1.3) \qquad \widehat{h}^F : \Omega(B; F) \to \Omega\left(B; \overline{F}^*\right).$$

Then the adjoint flat connection $\left(\nabla^F\right)^*$ on $F$ is given by

$$(1.4) \qquad \left(\nabla^F\right)^* = \left(\widehat{h}^F\right)^{-1} \nabla^{\overline{F}^*} \widehat{h}^F.$$

Define $\omega\left(\nabla^F, h^F\right) \in \Omega^1(B; \operatorname{End}(F))$ by

$$(1.5) \qquad \omega\left(\nabla^F, h^F\right) = \left(\nabla^F\right)^* - \nabla^F = \left(h^F\right)^{-1} \left(\nabla^F h^F\right).$$

In the rest of this section, except where otherwise indicated, we will abbreviate $\omega\left(\nabla^F, h^F\right)$ by $\omega$. With respect to a locally-defined covariantly-constant basis of $F$, $h^F$ is locally a Hermitian matrix-valued function on $B$ and we can write $\omega$ more simply as

$$(1.6) \qquad \omega = \left(h^F\right)^{-1} dh^F.$$

**Definition 1.1 :** The connection $\nabla^{F,u}$ on $F$ is given by

$$(1.7) \qquad \nabla^{F,u} = \nabla^F + \frac{\omega}{2}.$$

It is easy to see that $\nabla^{F,u}$ is compatible with $h^F$, and

$$(1.8) \qquad \nabla^{F,u} \omega = 0.$$

The curvature of $\nabla^{F,u}$ is given by

$$(1.9) \qquad \left(\nabla^{F,u}\right)^2 = -\frac{\omega^2}{4}.$$

**Definition 1.2 :** Define $P(\nabla^F, h^F) \in \Omega^{\text{even}}(B)$ by

$$(1.10) \qquad P(\nabla^F, h^F) = P\left(\frac{\omega^2}{8i\pi}\right).$$

**Lemma 1.3 :** *We have*

$$(1.11) \qquad P(\nabla^F, h^F) = P(0).$$

**Proof :** We can write $P(A)$ as a power series in the variables $\{n_j(A)\}_{j=1}^N$. It is easy to see algebraically that if $j > 0$ then $n_j(\nabla^F, h^F) = 0$. ∎

Let $z$ be an odd Grassmann variable, so that $z^2 = 0$. Given $\alpha \in \Omega(B) \widehat{\otimes} \mathbb{C}[z]$, we can write $\alpha$ in the form

(1.12) $$\alpha = \alpha_0 + z\alpha_1$$

with $\alpha_0, \alpha_1 \in \Omega(B)$. Put

(1.13) $$\text{Tr}_z[\alpha] = \alpha_1.$$

**Definition 1.4 :** Define $P^z(\nabla^F, h^F) \in \Omega^{\text{odd}}(B)$ by

(1.14) $$P^z(\nabla^F, h^F) = \text{Tr}_z\left[P(\frac{\omega^2}{8i\pi} + z\frac{\omega}{2})\right].$$

In particular,

(1.15) $$n_j^z(\nabla^F, h^F) = j \cdot 2^{-(2j-1)}(2i\pi)^{-(j-1)} \text{Tr}\left[\omega^{2j-1}\right].$$

**Lemma 1.5 :** *If $P$ and $Q$ are two ad-invariant power series then*

(1.16) $$(PQ)^z(\nabla^F, h^F) = P(0) \cdot Q^z(\nabla^F, h^F) + P^z(\nabla^F, h^F) \cdot Q(0).$$

**Proof :** This follows from Lemma 1.3. ∎

**Lemma 1.6 :** *The odd form $P^z(\nabla^F, h^F)$ is closed and its de Rham cohomology class is independent of $h^F$. If $P$ is real then $P^z(\nabla^F, h^F)$ is real.*

**Proof :** In the case $P(A) = n_j(A)$, $j > 0$, a simple proof of the statement of the lemma was given in [BLo, Theorems 1.8 and 1.11]. The general case follows from expressing $P(A)$ as a power series in the variables $\{n_j(A)\}_{j=1}^N$ and using Lemma 1.5. ∎

**Lemma 1.7 :** *We have*

(1.17) $$P^z(\nabla^{\overline{F}^*}, h^{\overline{F}^*}) = -P^z(\nabla^F, h^F).$$

**Proof :** One has

$$\omega\left(\nabla^F, h^F\right) = -\left(\widehat{h}^F\right)^{-1} \omega\left(\nabla^{\overline{F}^*}, h^{\overline{F}^*}\right) \widehat{h}^F. \tag{1.18}$$

The lemma follows from (1.14). ∎

**Definition 1.8 :** Let $P^z(\nabla^F) \in \mathrm{H}^{\mathrm{odd}}(B; \mathbb{R})$ denote the de Rham cohomology class of $P^z(\nabla^F, h^F)$.

**Remark 1.9 :** Given the polynomial $P$, there is a corresponding Cheeger-Chern-Simons class $\widetilde{P}(\nabla^F) \in \mathrm{H}^{\mathrm{odd}}(B; \mathbb{C}/\mathbb{Z})$ of the flat bundle $F$ [CS]. If $P(A) = n_j(A)$ then, up to a multiplicative constant, $P^z(\nabla^F)$ is the same as the imaginary part of $\widetilde{P}(\nabla^F)$ [BLo, Proposition 1.14].

**Theorem 1.10 :** *We have*

$$\omega^{2N} = 0. \tag{1.19}$$

*In particular, for $j > N$,*

$$n_j^z(\nabla^F, h^F) = 0. \tag{1.20}$$

*For $j \leq N$,*

$$c_j^z(\nabla^F, h^F) = \frac{(-1)^{j-1}}{j} \, n_j^z(\nabla^F, h^F). \tag{1.21}$$

**Proof :** Identity (1.19) appears in [R, Theorem 4.1]. We give a direct proof, which is essentially the same as that of [R]. By the Cayley-Hamilton theorem,

$$\left(\frac{\omega^2}{8i\pi}\right)^N + \sum_{j=1}^{N} (-1)^j c_j(\nabla^F, h^F) \cdot \left(\frac{\omega^2}{8i\pi}\right)^{N-j} = 0. \tag{1.22}$$

By Lemma 1.3, for $j > 0$,

$$c_j(\nabla^F, h^F) = 0. \tag{1.23}$$

Equation (1.19) follows.

Equation (1.20) follows from (1.15) and (1.19). Finally, if $j \le N$ then Newton's formula gives an identity of polynomials:

$$(1.24) \qquad n_j - c_1 n_{j-1} + \ldots + (-1)^{j-1} c_{j-1} n_1 + (-1)^j j c_j = 0.$$

Using Lemma 1.5, equation (1.21) follows. ∎

Let $E = E_+ \oplus E_-$ be a $\mathbb{Z}_2$-graded complex vector bundle on $B$. Let $\nabla^E = \nabla^{E_+} \oplus \nabla^{E_-}$ be a flat connection on $E$ which preserves the splitting $E = E_+ \oplus E_-$. Let $h^E = h^{E_+} \oplus h^{E_-}$ be a Hermitian metric on $E$ such that $E_+$ and $E_-$ are orthogonal. Put

$$(1.25) \qquad P^z(\nabla^E, h^E) = P^z(\nabla^{E_+}, h^{E_+}) - P^z(\nabla^{E_-}, h^{E_-}).$$

Given the flat vector bundle $F$, the associated flat vector bundle $\Lambda(\overline{F}^*)$ has a natural $\mathbb{Z}_2$-grading. If $h^F$ is a Hermitian metric on $F$ then there is an induced Hermitian metric $h^{\Lambda(\overline{F}^*)}$ on $\Lambda(\overline{F}^*)$.

**Theorem 1.11 :** *We have*

$$(1.26) \qquad \mathrm{ch}^z\left(\nabla^{\Lambda(\overline{F}^*)}, h^{\Lambda(\overline{F}^*)}\right) = \frac{1}{N} n_N^z(\nabla^F, h^F).$$

**Proof :** In general,

$$(1.27) \qquad \mathrm{ch}(\Lambda(A)) = \det(I - e^A) = (-1)^N \,\mathrm{Td}(-A)^{-1} \cdot \det(A).$$

Using (1.16),

$$(1.28) \qquad \mathrm{ch}^z\left(\nabla^{\Lambda(\overline{F}^*)}, h^{\Lambda(\overline{F}^*)}\right) = (-1)^N c_N^z\left(\nabla^{\overline{F}^*}, h^{\overline{F}^*}\right) = -\frac{1}{N} n_N^z\left(\nabla^{\overline{F}^*}, h^{\overline{F}^*}\right).$$

The theorem now follows from Lemma 1.7. ∎

**Corollary 1.12 :** *In* $\mathrm{H}^{\mathrm{odd}}(B; \mathbb{R})$, *one has the equality*

$$(1.29) \qquad \mathrm{ch}^z\left(\nabla^{\Lambda(\overline{F}^*)}\right) = \frac{1}{N} n_N^z(\nabla^F).$$

*In particular,* $\mathrm{ch}^z\left(\nabla^{\Lambda(\overline{F}^*)}\right)$ *is concentrated in degree* $2N - 1$.

**Proof :** This is an immediate consequence of Theorem 1.11. ∎

**c) Topological description of the characteristic classes**

The classes $n_j^z(\nabla^F)$ are the characteristic classes (of flat vector bundles) which are of interest to us. A more topological description of them can be given as follows. Let $V$ be a finite-dimensional complex vector space. Let $\mathrm{H}_c^*(GL(V); \mathbb{R})$ denote the continuous group cohomology of $GL(V)$, meaning the cohomology of the complex of Eilenberg-Maclane cochains on $GL(V)$ which are continous in their arguments. Let $GL(V)_\delta$ denote $GL(V)$ with the discrete topology and let $BGL(V)_\delta$ denote its classifying space. The cohomology group $\mathrm{H}^*(BGL(V)_\delta; \mathbb{R})$ is isomorphic to the (discrete) group cohomology $\mathrm{H}^*(GL(V); \mathbb{R})$. There is a forgetful map

$$(1.30) \qquad \mu_V : \mathrm{H}_c^*(GL(V); \mathbb{R}) \to \mathrm{H}^*(BGL(V)_\delta; \mathbb{R}).$$

Fix a basepoint $* \in B$. Put $\Gamma = \pi_1(B, *)$ and let $h : B \to B\Gamma$ be the classifying map for the universal cover of $B$, defined up to homotopy. Let $V$ be the fiber of $F$ above $*$. The holonomy of $F$ is a homomorphism $r : \Gamma \to GL(V)$, and induces a map $Br : B\Gamma \to BGL(V)_\delta$. Then the flat bundle $F$ is classified by the homotopy class of maps $\nu = Br \circ h : B \to BGL(V)_\delta$. One can show that there is a class $n_{j,V}^z \in \mathrm{H}^{2j-1}(GL(V); \mathbb{R})$ such that $n_j^z(\nabla^F) = \nu^*(n_{j,V}^z)$, and a class $N_{j,V}^z \in \mathrm{H}_c^{2j-1}(GL(V); \mathbb{R})$ such that $n_{j,V}^z = \mu_V(N_{j,V}^z)$. For example, $N_{1,V}^z$ is given by the homomorphism $g \to \ln|\det(g)|$ from $GL(V)$ to $(\mathbb{R}, +)$.

Put $G = GL(V)$ and $K = U(V)$. Denote the Lie algebras of $G$ and $K$ by $\gamma = gl(V)$ and $\kappa = u(V)$, respectively. The quotient space $\gamma/\kappa$ is isomorphic to the space of Hermitian endomorphisms of $V$, and carries an adjoint representation of $K$. One has that $\mathrm{H}_c^*(GL(V); \mathbb{R})$ is isomorphic to $\mathrm{H}^*(\gamma, K; \mathbb{R})$, the cohomology of the complex $\mathrm{C}^*(\gamma, K; \mathbb{R}) = \mathrm{Hom}_K(\Lambda^*(\gamma/\kappa), \mathbb{R})$ [BW, Chapter IX, §5]. In fact, the differential of this complex vanishes, and so $\mathrm{H}_c^*(GL(V); \mathbb{R}) = \mathrm{C}^*(\gamma, K; \mathbb{R})$ [BW, Chapter II, Corollary 3.2]. Thus the classes $\{n_j^z(\nabla^F)\}_{j=1}^\infty$ arise indirectly from $K$-invariant forms on $\gamma/\kappa$. It is possible to see the relationship between $n_j^z(\nabla^F)$ and $\mathrm{C}^{2j-1}(\gamma, K; \mathbb{R})$ more directly [BLo, §1g]. In particular, define a $(2j-1)$-form $\Phi_j$ on $\gamma/\kappa$ by sending Hermitian endomorphisms $M_1, \ldots, M_{2j-1}$ to

$$(1.31) \qquad \Phi_j(M_1, \ldots, M_{2j-1}) = \sum_{\sigma \in S_{2j-1}} (-1)^{\mathrm{sign}(\sigma)} \mathrm{Tr}\left[M_{\sigma(1)} \ldots M_{\sigma(2j-1)}\right].$$

Then $\Phi_j$ is an element of $\mathrm{C}^{2j-1}(\gamma, K; \mathbb{R})$ which, up to an overall multiplicative constant, corresponds to $N_{j,V}^z$.

The compact dual of the symmetric space $G/K$ is $G^d/K$, where $G^d = U(V) \times U(V)$. Let $\gamma^d = u(V) \oplus u(V)$ be the Lie algebra of $G^d$. Duality gives an isomorphism between $H^*(\gamma, K; \mathbb{R})$ and $H^*(\gamma^d, K; \mathbb{R}) = H^*(U(V); \mathbb{R}) = \Lambda(x_1, x_3, ..., x_{2\dim(V)-1})$. It follows that the classes $\{N^z_{j,V}\}_{j=1}^{\dim(V)}$ are algebraically independent.

If $V$ is the complexification of a real vector space $V_{\mathbb{R}}$ then one can apply the same arguments with $G = GL(V_{\mathbb{R}})$, $K = O(V_{\mathbb{R}})$ and $G^d = U(V)$. One obtains that if the flat complex vector bundle $F$ is the complexification of a flat real vector bundle $F_{\mathbb{R}}$ and $j$ is even then $n^z_j(\nabla^F)$ vanishes.

**II - Vanishing of a group cohomology class of $GL(N, \mathbb{Z})$**

In this section we use Berezin integrals to define certain forms $\delta_t$ and $\rho_t$ on the total space of a flat vector bundle, which are closed when pulled back to the base by a flat section. We also define transgressing forms $\epsilon_t$ and $\sigma_t$. We use these forms to prove Theorems 0.1 and 0.2.

The section is organized as follows. In a) we briefly review the Berezin integral. In b) we review the construction of Mathai-Quillen of the Thom form $\alpha_t$ of a real rank-$N$ vector bundle $\mathcal{V}$. If $\mathcal{V}$ is flat, we then define the forms $\delta_t \in \Omega^{2N-1}(\mathcal{V})$ and $\epsilon_t \in \Omega^{2N-2}(\mathcal{V})$, and establish their basic properties. In c) we construct useful auxiliary forms $\rho_t \in \Omega^{2N-1}(\mathcal{V})$ and $\sigma_t \in \Omega^{2N-2}(\mathcal{V})$. In d) we use these forms, along with the Poisson summation formula, to prove Theorems 0.1 and 0.2.

**a) Berezin integrals**

Let $V$ now be a real oriented inner-product space of dimension $N$. Let $\{\psi_k\}_{k=1}^N$ be an oriented orthonormal basis of $V$. We can identify $\Lambda(V)$ with $\Lambda(\psi_1, \ldots, \psi_N)$. The Berezin integral $\int^B : \Lambda(V) \to \mathbb{R}$ is defined to be the linear functional which vanishes on $\Lambda^k(V)$ unless $k = N$, in which case it is given by

$$(2.1) \qquad \int^B \psi_1 \psi_2 \ldots \psi_N = 1.$$

If $A$ is a graded-commutative superalgebra over $\mathbb{R}$, there is an extension of the Berezin integral to a linear map

$$(2.2) \qquad \int^B : A \,\widehat{\otimes}\, \Lambda(V) \to A$$

such that for $a \in A$ and $\alpha \in \Lambda(V)$,

$$(2.3) \qquad \int^B a \cdot \alpha = a \int^B \alpha.$$

Let $\widehat{V}$ be another copy of $V$. The exterior algebra

$$(2.4) \qquad \Lambda\left(V \oplus \widehat{V}\right) = \Lambda\left(\psi_1, \ldots, \psi_N, \widehat{\psi}_1, \ldots, \widehat{\psi}_N\right)$$

has a bigrading as

$$(2.5) \qquad \Lambda\left(V \oplus \widehat{V}\right) = \bigoplus_{k,l=1}^{n} \Lambda^k(V) \,\widehat{\otimes}\, \Lambda^l(\widehat{V}).$$

Then the Berezin integral $\int^B : \Lambda\left(V \oplus \widehat{V}\right) \to \mathbb{R}$ vanishes on $\Lambda^k(V) \widehat{\otimes} \Lambda^l(\widehat{V})$ unless $k = l = N$, in which case it is given by

$$\text{(2.6)} \qquad \int^B \psi_1 \ldots \psi_N \widehat{\psi}_1 \ldots \widehat{\psi}_N = 1.$$

The Berezin integral on $\Lambda\left(V \oplus \widehat{V}\right)$ is independent of the orientation of $V$.

**Lemma 2.1 :** *If $V$ is an antisymmetric $(N \times N)$-matrix with even entries and $W$ is a symmetric $(N \times N)$-matrix with odd entries then*

$$\text{(2.7)} \qquad \int^B \left\langle \psi, W\widehat{\psi} \right\rangle e^{\frac{1}{2}\langle \psi, V\psi \rangle - \frac{1}{2}\left\langle \widehat{\psi}, V\widehat{\psi} \right\rangle} = -\pi^{N-1} \text{Tr}_z \left[ \det\left(\frac{V}{i\pi} + zW\right) \right].$$

*If $N > 1$ then*

$$\text{(2.8)} \qquad \int^B \left\langle \psi, \widehat{\psi} \right\rangle e^{\frac{1}{2}\langle \psi, V\psi \rangle - \frac{1}{2}\left\langle \widehat{\psi}, V\widehat{\psi} \right\rangle} = 0$$

*and if $N = 1$ then*

$$\text{(2.9)} \qquad \int^B \left\langle \psi, \widehat{\psi} \right\rangle e^{\frac{1}{2}\langle \psi, V\psi \rangle - \frac{1}{2}\left\langle \widehat{\psi}, V\widehat{\psi} \right\rangle} = 1.$$

**Proof :** Equation (2.7) can be checked by putting $V$ into normal form. Equations (2.8) and (2.9) follow similarly. ∎

Let $A \in \text{End}(V)$ be antisymmetric. Then we can identify $A$ with an element of $\Lambda^2(V)$ by

$$\text{(2.10)} \qquad A \to \frac{1}{2} \left\langle \psi, A\psi \right\rangle.$$

**b) Thom-like forms**

Let $\mathcal{V}$ be a real rank-$N$ oriented vector bundle over a smooth connected manifold $B$, with projection map $\pi : \mathcal{V} \to B$. Let $h^\mathcal{V}$ be a metric on $\mathcal{V}$ and let $\nabla^\mathcal{V}$ be a compatible connection. The connection $\nabla^\mathcal{V}$ gives a splitting

$$\text{(2.11)} \qquad T\mathcal{V} = \pi^*TB \oplus \pi^*\mathcal{V}.$$

There is an induced connection $\nabla^{\Lambda(\mathcal{V})}$ on $\Lambda(\mathcal{V})$. Put $\mathcal{E} = \Lambda(\pi^*\mathcal{V})$, a $\mathbb{Z}$-graded vector bundle on $\mathcal{V}$ with connection $\nabla^\mathcal{E} = \pi^*\nabla^{\Lambda(\mathcal{V})}$. We first define the Thom-form of Mathai-Quillen [MQ], following the notation of [BGV, Section 1.6].

The Berezin integral gives a map $\int^B : \Omega(\mathcal{V}; \mathcal{E}) \to \Omega(\mathcal{V})$. There are certain elements of $\Omega(\mathcal{V}; \mathcal{E})$ of interest, namely

1. the tautological section $\mathbf{x} \in \Omega^0(\mathcal{V}; \Lambda^1(\pi^*\mathcal{V}))$;
2. the element $|\mathbf{x}|^2 \in \Omega^0(\mathcal{V}; \Lambda^0(\pi^*\mathcal{V}))$;
3. the element $\nabla^{\mathcal{E}} \mathbf{x} \in \Omega^1(\mathcal{V}; \Lambda^1(\pi^*\mathcal{V}))$;
4. the curvature $R^{\mathcal{V}} = \pi^* \left(\nabla^{\mathcal{V}}\right)^2$. Using (2.10), we can think of $R^{\mathcal{V}}$ as an element of $\Omega^2(\mathcal{V}; \pi^*\Lambda^2(\mathcal{V}))$.

Let $i_{\mathbf{x}} : \Omega^p(\mathcal{V}; \Lambda^q(\pi^*\mathcal{V})) \to \Omega^p(\mathcal{V}; \Lambda^{q-1}(\pi^*\mathcal{V}))$ be interior multiplication by the tautological section. For any $t > 0$ and $\alpha \in \Omega(\mathcal{V}; \mathcal{E})$, one has

$$(2.12) \qquad d \int^B \alpha = \int^B \left(\nabla^{\mathcal{E}} + 2\sqrt{t}\, i_{\mathbf{x}}\right) \alpha.$$

**Definition 2.2 :** For $t > 0$, define $A_t \in \Omega(\mathcal{V}; \mathcal{E})$ by

$$(2.13) \qquad A_t = \frac{R^{\mathcal{V}}}{2} + \sqrt{t}\, \nabla^{\mathcal{E}} \mathbf{x} + t\, |\mathbf{x}|^2.$$

**Definition 2.3 :** For $t > 0$, the Mathai-Quillen form $\alpha_t \in \Omega^N(\mathcal{V})$ is given by

$$(2.14) \qquad \alpha_t = (-1)^{\frac{N(N+1)}{2}} \pi^{-\frac{N}{2}} \int^B e^{-A_t}.$$

The form $\beta_t \in \Omega^{N-1}(\mathcal{V})$ is given by

$$(2.15) \qquad \beta_t = -(-1)^{\frac{N(N+1)}{2}} \pi^{-\frac{N}{2}} \int^B \frac{\mathbf{x}}{2\sqrt{t}} e^{-A_t}.$$

**Theorem 2.4 ([MQ, Thm. 7.6], [BGV, Prop. 1.53]) :**

a. The form $\alpha_t$ is closed.

b. In addition,

$$(2.16) \qquad \frac{\partial \alpha_t}{\partial t} = d\beta_t.$$

**Proof :** a. One can check that

$$(2.17) \qquad \left(\nabla^{\mathcal{E}} + 2\sqrt{t}\, i_{\mathbf{x}}\right) A_t = 0.$$

Using (2.12) and (2.17), part a. of the theorem follows.

Part b. of the theorem can be easily checked directly, but we will give a more general construction which will be of use later. Put $B' = B \times \mathbb{R}^+$ and $\mathcal{V}' = \mathcal{V} \times \mathbb{R}^+$. Define $\pi' : \mathcal{V}' \to B'$ by $\pi'(f, s) = (\pi(f), s)$. Let $\rho_B : B' \to B$ and $\rho_\mathcal{V} : \mathcal{V}' \to \mathcal{V}$ be the projection maps. Then $\mathcal{V}' = \rho_B^* \mathcal{V}$. Using the product structure on $\mathcal{V}'$, we can write exterior differentiation on $\Omega(\mathcal{V}')$ as

$$(2.18) \qquad d' = d + ds\, \partial_s.$$

Let $h^{\mathcal{V}'}$ be the metric on $\mathcal{V}'$ which restricts to $s \cdot h^\mathcal{V}$ on $\mathcal{V} \times \{s\}$. Then

$$(2.19) \qquad \nabla^{\mathcal{V}'} = \rho_B^* \nabla^\mathcal{V} + \frac{ds}{2s}$$

is a connection on $\mathcal{V}'$ which is compatible with $h^{\mathcal{V}'}$. Furthermore, $\left(\nabla^{\mathcal{V}'}\right)^2 = \rho_B^* \left(\nabla^\mathcal{V}\right)^2$.

We now apply the preceding formalism to the vector bundle $\mathcal{V}'$ with metric $h^{\mathcal{V}'}$ and connection $\nabla^{\mathcal{V}'}$. Using an obvious notation, $A'_1 \in \Omega(\mathcal{V}'; \mathcal{E}')$ is given by

$$(2.20) \qquad A'_1 = \frac{R^\mathcal{V}}{2} + \sqrt{s}\left(\nabla^\mathcal{E} \mathbf{x} + ds \wedge \frac{\mathbf{x}}{2s}\right) + s|\mathbf{x}|^2 = A_s + ds \wedge \frac{\mathbf{x}}{2\sqrt{s}}.$$

Then

$$(2.21) \qquad \alpha'_1 = \int^B e^{-A'_1} = \left(\int^B e^{-A_s}\right) - ds \wedge \int^B \frac{\mathbf{x}}{2\sqrt{s}} e^{-A_s}.$$

From part a. of the theorem,

$$(2.22) \qquad d'\alpha'_1 = 0.$$

Part b. of the theorem now follows from (2.18), (2.21) and (2.22). ∎

**Remark 2.5 :** The closed form $\alpha_t$ is a Thom form in the sense that it is rapidly decreasing at infinity and the fiberwise integral $\pi_* \alpha_t$ is identically 1 on $B$.

We now assume that $\mathcal{V}$ has a flat connection $\nabla^\mathcal{V}$. Let $h^\mathcal{V}$ be a metric on $\mathcal{V}$. Define $\omega \in \Omega^1(B; \text{End}(\mathcal{V}))$ as in (1.5) and $\nabla^{\mathcal{V},u}$ as in (1.7).

**Lemma 2.6 :** *The $j$-form $\omega^j$ takes value in the symmetric endomorphisms of $\mathcal{V}$ if $j \equiv 0, 1 \pmod{4}$ and in the antisymmetric endomorphisms of $\mathcal{V}$ if $j \equiv 2, 3 \pmod{4}$.*

**Proof:** One sees from (1.5) that $\omega$ takes value in the symmetric endomorphisms of $\mathcal{V}$. Then $\left(\omega^j\right)^T = (-1)^{\frac{j(j-1)}{2}} \omega^j$. ∎

Let $\widehat{\mathcal{V}}$ be another copy of $\mathcal{V}$. Put $\widehat{\mathcal{E}} = \Lambda\left(\pi^*\widehat{\mathcal{V}}\right)$. Following (2.6), the Berezin integral gives a map $\int^B : \Omega\left(\mathcal{V}; \mathcal{E}\widehat{\otimes}\widehat{\mathcal{E}}\right) \to \Omega\left(\mathcal{V}\right)$.

There are certain elements of $\Omega\left(\mathcal{V}; \mathcal{E}\widehat{\otimes}\widehat{\mathcal{E}}\right)$ of interest, namely
1. the tautological section $\mathbf{x} \in \Omega^0(\mathcal{V}; \Lambda^1(\pi^*\mathcal{V}) \widehat{\otimes} \Lambda^0(\pi^*\widehat{\mathcal{V}}))$;
2. the tautological section $\widehat{\mathbf{x}} \in \Omega^0(\mathcal{V}; \Lambda^0(\pi^*\mathcal{V}) \widehat{\otimes} \Lambda^1(\pi^*\widehat{\mathcal{V}}))$;
3. the element $|\mathbf{x}|^2 = |\widehat{\mathbf{x}}|^2 \in \Omega^0(\mathcal{V}; \Lambda^0(\pi^*\mathcal{V}) \widehat{\otimes} \Lambda^0(\pi^*\widehat{\mathcal{V}}))$;
4. the element $\left\langle \psi, \widehat{\psi} \right\rangle \in \Omega^0(\mathcal{V}; \Lambda^1(\pi^*\mathcal{V}) \widehat{\otimes} \Lambda^1(\pi^*\widehat{\mathcal{V}}))$;
5. the element $\nabla^{\mathcal{E}\widehat{\otimes}\widehat{\mathcal{E}},u}\mathbf{x} \in \Omega^1(\mathcal{V}; \Lambda^1(\pi^*\mathcal{V}) \widehat{\otimes} \Lambda^0(\pi^*\widehat{\mathcal{V}}))$;
6. the element $\nabla^{\mathcal{E}\widehat{\otimes}\widehat{\mathcal{E}},u}\widehat{\mathbf{x}} \in \Omega^1(\mathcal{V}; \Lambda^0(\pi^*\mathcal{V}) \widehat{\otimes} \Lambda^1(\pi^*\widehat{\mathcal{V}}))$;
7. the curvature $R^{\mathcal{V},u} = \left(\pi^*\nabla^{\mathcal{V},u}\right)^2 = -\frac{\pi^*\omega^2}{4}$. We will think of $R^{\mathcal{V},u}$ as an element of $\Omega^2(\mathcal{V}; \Lambda^2(\pi^*\mathcal{V}) \widehat{\otimes} \Lambda^0(\pi^*\widehat{\mathcal{V}}))$, namely $R^{\mathcal{V},u} = -\frac{1}{8}\left\langle \psi, (\pi^*\omega^2)\psi \right\rangle$.
8. the element $\widehat{\omega} \in \Omega^1(\mathcal{V}; \Lambda^1(\pi^*\mathcal{V}) \widehat{\otimes} \Lambda^1(\pi^*\widehat{\mathcal{V}}))$ given by $\widehat{\omega} = \left\langle \psi, (\pi^*\omega)\widehat{\psi} \right\rangle$;
9. the element $\widehat{\omega^2} \in \Omega^2(\mathcal{V}; \Lambda^0(\pi^*\mathcal{V}) \widehat{\otimes} \Lambda^2(\pi^*\widehat{\mathcal{V}}))$ given by $\widehat{\omega^2} = \frac{1}{2}\left\langle \widehat{\psi}, (\pi^*\omega^2)\widehat{\psi} \right\rangle$.

As in (2.12), for any $t > 0$ and $\alpha \in \Omega\left(\mathcal{V}; \mathcal{E} \widehat{\otimes} \widehat{\mathcal{E}}\right)$, one has

$$\text{(2.23)} \qquad d\int^B \alpha = \int^B \left(\nabla^{\mathcal{E}\widehat{\otimes}\widehat{\mathcal{E}},u} + 2\sqrt{t}\, i_{\mathbf{x}}\right)\alpha.$$

**Definition 2.7:** For $t > 0$, define $B_t \in \Omega\left(\mathcal{V}; \mathcal{E} \widehat{\otimes} \widehat{\mathcal{E}}\right)$ by

$$\text{(2.24)} \qquad B_t = \frac{R^{\mathcal{V},u}}{2} + \sqrt{t}\, \nabla^{\mathcal{E}\widehat{\otimes}\widehat{\mathcal{E}},u}\mathbf{x} + t\,|\mathbf{x}|^2 + \frac{1}{8}\widehat{\omega^2}.$$

**Definition 2.8:** For $t > 0$, define $\delta_t \in \Omega^{2N-1}(\mathcal{V})$ by

$$\text{(2.25)} \qquad \delta_t = \int^B \left(\frac{1}{4}\widehat{\omega} - \sqrt{t}\,\widehat{\mathbf{x}}\right) e^{-B_t}$$

and $\epsilon_t \in \Omega^{2N-2}(\mathcal{V})$ by

$$\text{(2.26)} \qquad \epsilon_t = -\int^B \left(\frac{1}{4t}\left\langle \psi, \widehat{\psi}\right\rangle + \frac{\mathbf{x}}{2\sqrt{t}}\left(\frac{1}{4}\widehat{\omega} - \sqrt{t}\,\widehat{\mathbf{x}}\right)\right) e^{-B_t}.$$

**Theorem 2.9 :**

a. Let $U$ be an open subset of $B$ and let $\lambda : U \to \mathcal{V}$ be a flat section of $\mathcal{V}$ over $U$. Then for $t > 0$,

$$d\left(\lambda^* \delta_t\right) = 0. \tag{2.27}$$

b. In addition,

$$\frac{\partial \left(\lambda^* \delta_t\right)}{\partial t} = d\left(\lambda^* \epsilon_t\right). \tag{2.28}$$

**Proof :** a. As in (2.17), we have

$$\left(\nabla^{\mathcal{E}\widehat{\otimes}\widehat{\mathcal{E}},u} + 2\sqrt{t}\, i_{\mathbf{x}}\right)\left(\frac{R^{\mathcal{V},u}}{2} + \sqrt{t}\, \nabla^{\mathcal{E}\widehat{\otimes}\widehat{\mathcal{E}},u}\mathbf{x} + t\, |\mathbf{x}|^2\right) = 0. \tag{2.29}$$

From (1.8), we have

$$\left(\nabla^{\mathcal{E}\widehat{\otimes}\widehat{\mathcal{E}},u} + 2\sqrt{t}\, i_{\mathbf{x}}\right)\widehat{\omega^2} = 0. \tag{2.30}$$

Furthermore,

$$\left(\nabla^{\mathcal{E}\widehat{\otimes}\widehat{\mathcal{E}},u} + 2\sqrt{t}\, i_{\mathbf{x}}\right)\left(\frac{1}{4}\widehat{\omega} - \sqrt{t}\,\widehat{\mathbf{x}}\right) = -\sqrt{t}\, \nabla^{\mathcal{E}\widehat{\otimes}\widehat{\mathcal{E}},u}\widehat{\mathbf{x}} + \frac{\sqrt{t}}{2}\, i_{\mathbf{x}}\widehat{\omega}. \tag{2.31}$$

Using (2.29), (2.30) and (2.31), we obtain

$$d\delta_t = \int^B \left(-\sqrt{t}\, \nabla^{\mathcal{E}\widehat{\otimes}\widehat{\mathcal{E}},u}\widehat{\mathbf{x}} + \frac{\sqrt{t}}{2}\, i_{\mathbf{x}}\widehat{\omega}\right) e^{-B_t}. \tag{2.32}$$

As elements of $\Omega^1(U; \Lambda^0(\mathcal{V}) \widehat{\otimes} \Lambda^1(\widehat{\mathcal{V}}))$, there is an equality

$$\lambda^* \left(\nabla^{\mathcal{E}\widehat{\otimes}\widehat{\mathcal{E}},u}\widehat{\mathbf{x}}\right) = \frac{1}{2}\, \lambda^* \left(i_{\mathbf{x}}\widehat{\omega}\right). \tag{2.33}$$

Equation (2.27) follows.

To prove b., we continue with the setup of the proof of Theorem 2.4.b. Let $\nabla^{\mathcal{V}'}$ now be the flat connection $\rho_B^* \nabla^{\mathcal{V}}$. Then $\lambda' = \rho_B^* \lambda$ is a flat section. We have

$$\left(\nabla^{\mathcal{V}'}\right)^* = \rho_B^* \left(\nabla^{\mathcal{V}}\right)^* + \frac{ds}{s}. \tag{2.34}$$

and

$$\omega' = \rho_B^* \omega + \frac{ds}{s}. \tag{2.35}$$

Then

(2.36)
$$(\omega')^2 = \rho_B^* \omega^2$$

and

(2.37)
$$\widehat{\omega}' = \rho_{\mathcal{V}}^* \widehat{\omega} + \left\langle \psi, \frac{ds}{s} \widehat{\psi} \right\rangle = \rho_{\mathcal{V}}^* \widehat{\omega} - \frac{ds}{s} \left\langle \psi, \widehat{\psi} \right\rangle.$$

The unitary connection on $\mathcal{V}'$ is given by

(2.38)
$$\nabla^{\mathcal{V}',u} = \rho_B^* \nabla^{\mathcal{V},u} + \frac{ds}{2s}.$$

Using an obvious notation, $B_1' \in \Omega\left(\mathcal{V}'; \mathcal{E}' \widehat{\otimes} \widehat{\mathcal{E}'}\right)$ is then given by

(2.39)
$$B_1' = B_s + ds \wedge \frac{\mathbf{x}}{2\sqrt{s}}.$$

Thus

(2.40)
$$\delta_1' = \int^B \left(\frac{1}{4} \widehat{\omega}' - \sqrt{s} \, \widehat{\mathbf{x}}\right) e^{-B_1'} = \int^B \left(\frac{1}{4} \widehat{\omega} - \sqrt{s} \, \widehat{\mathbf{x}}\right) e^{-B_s}$$
$$- ds \wedge \int^B \left(\frac{1}{4s} \left\langle \psi, \widehat{\psi} \right\rangle + \frac{\mathbf{x}}{2\sqrt{s}} \left(\frac{1}{4} \widehat{\omega} - \sqrt{s} \, \widehat{\mathbf{x}}\right)\right) e^{-B_s}.$$

From part a. of the theorem,

(2.41)
$$d'\left(\lambda'^* \delta_1'\right) = 0.$$

Part b. of the theorem now follows from (2.18), (2.40) and (2.41). ∎

**Remark 2.10 :** If $N$ is odd then $\alpha_t = \beta_t = 0$. If $N$ is even then $\delta_t = \epsilon_t = 0$.

**Remark 2.11 :** There is an evident analogy between the properties of $(\alpha_t, \beta_t)$ and $(\delta_t, \epsilon_t)$. However, there is the important difference that $\alpha_t$ can be pulled-back by an arbitrary section of $\mathcal{V}$ to get a closed form on $B$, whereas $\delta_t$ can only be pulled-back by a flat section of $\mathcal{V}$ if one wants to get a closed form on $B$.

**c) Some auxiliary forms**

We use the notation of Section 2b. Suppose that the holonomy of the flat connection $\nabla^{\mathcal{V}}$ preserves a volume form $\eta$ on the fibers. We will assume that the metric $h^{\mathcal{V}}$ is such that the volume form on the fibers induced by $h^{\mathcal{V}}$ equals $\eta$. The holonomy of the adjoint flat connection $(\nabla^{\mathcal{V}})^*$ also preserves $\eta$.

**Definition 2.12 :** Define $\mathrm{Vol} \in \Omega^N(\mathcal{V})$ to be the extension of $\eta$ to $\mathcal{V}$, using $(\nabla^{\mathcal{V}})^*$.

More precisely, if $U$ is a contractible open set in $B$ and $U \times \mathbb{R}^N$ is a trivialization of $\mathcal{V}$ which is covariantly-constant with respect to $(\nabla^{\mathcal{V}})^*$, let $\phi : U \times \mathbb{R}^N \to \mathbb{R}^N$ be projection onto the fiber. Then on $U \times \mathbb{R}^N$, one has $\mathrm{Vol} = \phi^* \eta$. Clearly Vol is closed. Consider the Berezin integral $\int^B : \Omega(\mathcal{V}; \mathcal{E}) \to \Omega(\mathcal{V})$.

**Theorem 2.13 :** *We have*

$$(2.42) \qquad \mathrm{Vol} = (-1)^{\frac{N(N-1)}{2}} \int^B e^{(\nabla^{\mathcal{E}})^* \mathbf{x}}.$$

*If $\lambda : U \to \mathcal{V}$ is a section of $\mathcal{V}$ which is flat with respect to $\nabla^{\mathcal{V}}$ then*

$$(2.43) \qquad \lambda^* \mathrm{Vol} = (\omega \lambda)_1 \wedge \ldots \wedge (\omega \lambda)_N.$$

**Proof :** In terms of the trivialization $U \times \mathbb{R}^N$ above, we can write $(\nabla^{\mathcal{E}})^* \mathbf{x} = \langle d\phi, \psi \rangle$. Then

$$(2.44) \qquad (-1)^{\frac{N(N-1)}{2}} \int^B e^{\langle d\phi, \psi \rangle} = (-1)^{\frac{N(N-1)}{2}} \int^B \prod_{j=1}^N d\phi_j \cdot \psi_j = d\phi_1 \wedge \ldots \wedge d\phi_N.$$

Equation (2.42) follows. As elements of $\Omega^1(U; \mathcal{V})$, we have

$$(2.45) \qquad \lambda^*((\nabla^{\mathcal{E}})^* \mathbf{x}) = (\nabla^{\mathcal{E}} + \omega) \lambda = \omega \lambda.$$

Equation (2.43) follows. ∎

Now consider the Berezin integral $\int^B : \Omega\left(\mathcal{V}; \widehat{\mathcal{E}}\right) \to \Omega(\mathcal{V})$.

**Definition 2.14 :** For $t > 0$, define $\rho_t \in \Omega^{2N-1}(\mathcal{V})$ by

$$(2.46) \qquad \rho_t = e^{-\frac{|\widehat{\mathbf{x}}|^2}{4t}} t^{-N} \, \mathrm{Vol} \cdot \int^B \frac{\widehat{\mathbf{x}}}{\sqrt{t}} e^{-\frac{1}{8}\widehat{\omega}^2}$$

and $\sigma_t \in \Omega^{2N-2}(\mathcal{V})$ by

$$(2.47) \qquad \sigma_t = -e^{-\frac{|\widehat{\mathbf{x}}|^2}{4t}} t^{-N-1} \, i_x \mathrm{Vol} \cdot \int^B \frac{\widehat{\mathbf{x}}}{\sqrt{t}} e^{-\frac{1}{8}\widehat{\omega}^2}.$$

**Theorem 2.15 :**

a. Let $U$ be an open subset of $B$ and let $\lambda : U \to \mathcal{V}$ be a flat section of $\mathcal{V}$ over $U$. Then for $t > 0$,

(2.48) $$d(\lambda^* \rho_t) = 0.$$

b. In addition,

(2.49) $$\frac{\partial(\lambda^* \rho_t)}{\partial t} = d(\lambda^* \sigma_t).$$

**Proof :** If $N$ is even then $\rho_t = \sigma_t = 0$ and the theorem is trivially true. Thus we may assume that $N$ is odd.

a. We have

(2.50) $$d\rho_t = -e^{-\frac{|\widehat{\mathbf{x}}|^2}{4t}} t^{-N} \operatorname{Vol} \cdot \left( \frac{i_{\widehat{\mathbf{x}}}(\nabla^{\widehat{\mathcal{E}},u}\widehat{\mathbf{x}})}{2t} \int^B \frac{\widehat{\mathbf{x}}}{\sqrt{t}} e^{-\frac{1}{8}\widehat{\omega}^2} + \int^B \frac{\nabla^{\widehat{\mathcal{E}},u}\widehat{\mathbf{x}}}{\sqrt{t}} e^{-\frac{1}{8}\widehat{\omega}^2} \right).$$

Furthermore,

(2.51) $$\lambda^*\left(\nabla^{\widehat{\mathcal{E}},u}\widehat{\mathbf{x}}\right) = \nabla^{\mathcal{V},u}\lambda = \left(\nabla^{\mathcal{V}} + \frac{\omega}{2}\right)\lambda = \frac{\omega}{2}\lambda$$

and

(2.52) $$\lambda^*\left(i_{\widehat{\mathbf{x}}}(\nabla^{\widehat{\mathcal{E}},u}\widehat{\mathbf{x}})\right) = -\frac{1}{2}\langle \lambda, \omega\lambda \rangle.$$

Equation (2.48) now follows from the explicit representation of $\lambda^*\operatorname{Vol}$ as an $N$-form in (2.43).

b. We continue with the notation of the proof of Theorem 2.4.b, except that we now take $h^{\mathcal{V}'}$ to be the metric on $\mathcal{V}'$ which restricts to $\frac{1}{s} \cdot h^{\mathcal{V}}$ on $\mathcal{V} \times \{s\}$. Let $\nabla^{\mathcal{V}'}$ be the flat connection $\rho_B^* \nabla^{\mathcal{V}}$. Then

(2.53) $$\left(\nabla^{\mathcal{V}'}\right)^* = \rho_B^*\left(\nabla^{\mathcal{V}}\right)^* - \frac{ds}{s}.$$

Using (2.42) and (2.53), we have

(2.54) $$\operatorname{Vol}' = (-1)^{\frac{N(N-1)}{2}} \int^B e^{(\nabla^{\mathcal{V}})^*\mathbf{x} - \frac{ds}{s}\mathbf{x}} = \operatorname{Vol} - \frac{ds}{s} \wedge i_x\operatorname{Vol}.$$

Then

(2.55) $$\rho_1' = e^{-\frac{|\widehat{\mathbf{x}}|^2}{4s}} s^{-N} \operatorname{Vol} \cdot \int^B \frac{\widehat{\mathbf{x}}}{\sqrt{s}} e^{-\frac{1}{8}\widehat{\omega}^2} - ds \wedge e^{-\frac{|\widehat{\mathbf{x}}|^2}{4s}} s^{-N-1} i_x\operatorname{Vol} \cdot \int^B \frac{\widehat{\mathbf{x}}}{\sqrt{s}} e^{-\frac{1}{8}\widehat{\omega}^2}.$$

From part a. of the theorem,

$$(2.56) \qquad d'\left({\lambda'}^* \rho'_1\right) = 0.$$

Part b. of the theorem now follows from (2.18), (2.55) and (2.56). ∎

### d) Proof of the main theorem

Let $N$ be a positive odd integer. Let $B$ be a smooth connected manifold. Let $\widetilde{B}$ be the universal cover of $B$ and put $\Gamma = \pi_1(B)$. Let $\rho : \Gamma \to SL(N, \mathbb{Z})$ be a homomorphism. When convenient, we will consider $SL(N, \mathbb{Z})$ to be a subgroup of $SL(N, \mathbb{R})$ or $SL(N, \mathbb{C})$. Define a flat real rank-$N$ vector bundle on $B$:

$$(2.57) \qquad E = \widetilde{B} \times_\rho \mathbb{R}^N.$$

Let $\Lambda \subset E$ be the lattice

$$(2.58) \qquad \Lambda = \widetilde{B} \times_\rho \mathbb{Z}^N.$$

Let $E^*$ be the dual vector bundle to $E$ and let $\Lambda^* \subset E^*$ be the dual lattice

$$(2.59) \qquad \Lambda^* = \{\mu \in E^* : \forall \lambda \in \Lambda, \langle \mu, \lambda \rangle \in 2\pi \mathbb{Z}\}.$$

Let $\nabla^E$ be the canonical flat connection on $E$. Then $\nabla^E$ preserves the lattice $\Lambda$. The dual flat connection $\nabla^{E^*}$ preserves $\Lambda^*$.

From (2.57), there are volume forms on the fibers of $E$ coming from the standard volume form on $\mathbb{R}^N$. Let $\mathrm{Vol}(E/\Lambda)$ denote the common volume of the quotients of the fibers by $\Lambda$. Choose an inner product $h^E$ on $E$ which is compatible with these volume forms. We write $\omega_E$ for $\omega\left(\nabla^E, h^E\right) \in \Omega^1(B; \mathrm{End}(E))$ and $\omega_{E^*}$ for $\omega\left(\nabla^{E^*}, h^{E^*}\right) \in \Omega^1(B; \mathrm{End}(E^*))$

We now apply the formalism of Section 2b to the case $\mathcal{V} = E^*$. Define $\delta_t \in \Omega^{2N-1}(E^*)$ as in (2.25) and $\epsilon_t \in \Omega^{2N-2}(E^*)$ as in (2.26). If $U$ is a contractible open subset of $B$ then over $U$, the lattice $\Lambda$ consists of a countable number of disjoint copies of $U$. Thus the forms $\sum_{\mu \in \Lambda^*} \mu^* \delta_t$ and $\sum_{\mu \in \Lambda^*} \mu^* \epsilon_t$ are well-defined on $U$, and we obtain global forms $\sum_{\mu \in \Lambda^*} \mu^* \delta_t \in \Omega^{2N-1}(B)$ and $\sum_{\mu \in \Lambda^*} \mu^* \epsilon_t \in \Omega^{2N-2}(B)$.

**Theorem 2.16 :** *We have*

$$(2.60) \qquad d \sum_{\mu \in \Lambda^*} \mu^* \delta_t = 0$$

*and*

$$\frac{\partial}{\partial t} \sum_{\mu \in \Lambda^*} \mu^* \delta_t = d \sum_{\mu \in \Lambda^*} \mu^* \epsilon_t. \tag{2.61}$$

**Proof :** This is a consequence of Theorem 2.9. ∎

**Theorem 2.17 :** *If $K$ is a compact subset of $B$ then there is a constant $c > 0$ such that on $K$, as $t \to \infty$,*

$$\sum_{\mu \in \Lambda^*} \mu^* \delta_t = \frac{1}{2} \pi^{N-1} c_N^z \left( \nabla^E, h^E \right) + \mathcal{O}\left(e^{-ct}\right). \tag{2.62}$$

*If $N > 1$ then*

$$\sum_{\mu \in \Lambda^*} \mu^* \epsilon_t = \mathcal{O}\left(e^{-ct}\right) \tag{2.63}$$

*and if $N = 1$ then*

$$\sum_{\mu \in \Lambda^*} \mu^* \epsilon_t = -\frac{1}{4t} + \mathcal{O}\left(e^{-ct}\right). \tag{2.64}$$

**Proof :** It is enough to consider only the contribution of $\mu = 0$, as the other terms will be exponentially damped in $t$. From (2.25),

$$0^* \delta_t = \frac{1}{4} \int^B \left\langle \psi, \omega_{E^*} \widehat{\psi} \right\rangle e^{\frac{1}{16}\left\langle \psi, \omega_{E^*}^2 \psi \right\rangle - \frac{1}{16}\left\langle \widehat{\psi}, \omega_{E^*}^2 \widehat{\psi} \right\rangle}. \tag{2.65}$$

Using Lemma 2.1 with $V = \frac{\omega_{E^*}^2}{8}$ and $W = \frac{\omega_{E^*}}{2}$, we obtain

$$0^* \delta_t = -\frac{1}{2} \pi^{N-1} c_N^z \left( \nabla^{E^*}, h^{E^*} \right). \tag{2.66}$$

Equation (2.62) now follows from Lemma 1.7. Equations (2.63) and (2.64) follow similarly. ∎

Applying the results of Section 2c to the case $\mathcal{V} = E$, define $\rho_t \in \Omega^{2N-1}(E)$ as in (2.46) and $\sigma_t \in \Omega^{2N-2}(E)$ as in (2.47).

**Theorem 2.18 :** *The forms $\sum_{m \in \Lambda} m^* \rho_t \in \Omega^{2N-1}(B)$ and $\sum_{m \in \Lambda} m^* \sigma_t \in \Omega^{2N-2}(B)$ satisfy*

$$d \sum_{m \in \Lambda} m^* \rho_t = 0 \tag{2.67}$$

*and*

$$\frac{\partial}{\partial t} \sum_{m \in \Lambda} m^* \rho_t = d \sum_{m \in \Lambda} m^* \sigma_t. \tag{2.68}$$

**Proof :** This is a consequence of (2.48) and (2.49). ∎

**Theorem 2.19 :** *We have*

$$\sum_{\mu \in \Lambda^*} \mu^* \delta_t = 2^{-3N-1} \, \pi^{-\frac{N}{2}} \, \text{Vol}(E/\Lambda) \sum_{m \in \Lambda} m^* \rho_t \tag{2.69}$$

*and*

$$\sum_{\mu \in \Lambda^*} \mu^* \epsilon_t = 2^{-3N-1} \, \pi^{-\frac{N}{2}} \, \text{Vol}(E/\Lambda) \sum_{m \in \Lambda} m^* \sigma_t. \tag{2.70}$$

**Proof :** From (2.24), we have

$$\begin{aligned} \mu^* B_t &= -\frac{1}{16} \langle \psi, \omega_{E^*}^2 \psi \rangle + \sqrt{t} \left\langle \frac{1}{2} \omega_{E^*} \mu, \psi \right\rangle + t \, |\mu|^2 + \frac{1}{16} \langle \widehat{\psi}, \omega_{E^*}^2 \widehat{\psi} \rangle \\ &= \frac{1}{16} |\omega_{E^*} \psi|^2 + \sqrt{t} \left\langle \frac{1}{2} \mu, \omega_{E^*} \psi \right\rangle + t \, |\mu|^2 + \frac{1}{16} \langle \widehat{\psi}, \omega_{E^*}^2 \widehat{\psi} \rangle \\ &= t \left| \mu + \frac{1}{4\sqrt{t}} \omega_{E^*} \psi \right|^2 + \frac{1}{16} \langle \widehat{\psi}, \omega_{E^*}^2 \widehat{\psi} \rangle. \end{aligned} \tag{2.71}$$

Furthermore,

$$\mu^* \left( \frac{1}{4} \widehat{\omega} - \sqrt{t} \, \widehat{\mathbf{x}} \right) = -\sqrt{t} \left\langle \left( \mu + \frac{1}{4\sqrt{t}} \omega_{E^*} \psi \right), \widehat{\psi} \right\rangle. \tag{2.72}$$

Let $z$ be an auxiliary odd variable, satisfying $z^2 = 0$. Define $\text{Tr}_z$ as in (1.13). Then from equation (2.71), in terms of the Berezin integral $\int^B : \Omega\left(B; \Lambda(E^*) \widehat{\otimes} \Lambda(\widehat{E^*})\right) \to \Omega(B)$, we have

$$\mu^* \delta_t = \text{Tr}_z \left[ \int^B e^{-t \left| \mu + \frac{1}{4\sqrt{t}} \omega_{E^*} \psi + \frac{1}{2\sqrt{t}} z \widehat{\psi} \right|^2 - \frac{1}{16} \langle \widehat{\psi}, \omega_{E^*}^2 \widehat{\psi} \rangle } \right]. \tag{2.73}$$

From the Poisson summation formula, one has in general that if $b \in E^*$ then

$$\sum_{\mu \in \Lambda^*} e^{-t \, |\mu+b|^2} = (4\pi t)^{-\frac{N}{2}} \, \text{Vol}(E/\Lambda) \sum_{m \in \Lambda} e^{-\frac{|m|^2}{4t}} e^{i \langle m, b \rangle}. \tag{2.74}$$

Applying this to (2.73), we obtain

$$(2.75) \quad \sum_{\mu \in \Lambda^*} \mu^* \delta_t = 2^{-N} \pi^{-\frac{N}{2}} t^{-\frac{N}{2}} \operatorname{Vol}(E/\Lambda) \cdot$$

$$\operatorname{Tr}_z \left[ \int^B \sum_{m \in \Lambda} e^{-\frac{|m|^2}{4t}} e^{i\langle m, \frac{1}{4\sqrt{t}} \omega_{E^*} \psi + \frac{1}{2\sqrt{t}} z\widehat{\psi}\rangle - \frac{1}{16}\langle \widehat{\psi}, \omega_{E^*}^2 \widehat{\psi}\rangle} \right]$$

$$= 2^{-N-1} \pi^{-\frac{N}{2}} t^{-\frac{N}{2}} i \operatorname{Vol}(E/\Lambda) \cdot$$

$$\sum_{m \in \Lambda} e^{-\frac{|m|^2}{4t}} \int^B e^{i\langle m, \frac{1}{4\sqrt{t}} \omega_{E^*} \psi\rangle} \frac{\widehat{m}}{\sqrt{t}} e^{-\frac{1}{16}\langle \widehat{\psi}, \omega_{E^*}^2 \widehat{\psi}\rangle}$$

Define new Grassmann variables by $\eta = \left(\widehat{h}^E\right)^{-1} \psi$ and $\widehat{\eta} = \left(\widehat{h}^E\right)^{-1} \widehat{\psi}$. Using (1.18), we have

$$(2.76) \quad \langle m, \omega_{E^*} \psi \rangle = -\langle m, \omega_E \eta \rangle$$

and

$$(2.77) \quad \left\langle \widehat{\psi}, \omega_{E^*}^2 \widehat{\psi} \right\rangle = \langle \widehat{\eta}, \omega_E^2 \widehat{\eta} \rangle.$$

From (2.75), (2.76) and (2.77), in terms of the Berezin integral $\int^B : \Omega\left(B; \Lambda(E) \widehat{\otimes} \Lambda(\widehat{E})\right) \to \Omega(B)$, we have

$$(2.78) \quad \sum_{\mu \in \Lambda^*} \mu^* \delta_t = 2^{-3N-1} \pi^{-\frac{N}{2}} \operatorname{Vol}(E/\Lambda) \cdot$$

$$\sum_{m \in \Lambda} e^{-\frac{|m|^2}{4t}} t^{-N} \int^B (-1)^{\frac{N(N-1)}{2}} e^{\langle m, \omega_E \eta\rangle} \frac{\widehat{m}}{\sqrt{t}} e^{-\frac{1}{16}\langle \widehat{\eta}, \omega_E^2 \widehat{\eta}\rangle}$$

$$= 2^{-3N-1} \pi^{-\frac{N}{2}} \operatorname{Vol}(E/\Lambda) \sum_{m \in \Lambda} m^* \rho_t.$$

This proves (2.69). In view of the constructions of the proofs of Theorems 2.9.b and 2.15.b, equation (2.70) follows automatically from (2.69). ∎

**Corollary 2.20 :** *If $K$ is a compact subset of $B$ then there is a constant $c > 0$ such that on $K$, as $t \to 0$,*

$$(2.79) \quad \sum_{\mu \in \Lambda^*} \mu^* \delta_t = \mathcal{O}\left(e^{-\frac{c}{t}}\right)$$

*and*

$$\sum_{\mu \in \Lambda^*} \mu^* \epsilon_t = \mathcal{O}\left(e^{-\frac{c}{t}}\right). \tag{2.80}$$

**Proof :** As the sums in (2.69) and (2.70) can be taken over nonzero $m$, the theorem follows. ∎

**Theorem 2.21** *Let $B$ be a smooth connected manifold. Let $N$ be a positive odd integer and let $E$ be a flat real rank-$N$ vector bundle over $B$ with structure group contained in $GL(N, \mathbb{Z})$. Then $n_N^z(\nabla^E)$ vanishes in $\mathrm{H}^{2N-1}(B; \mathbb{R})$.*

**Proof :** Assume first that the structure group is contained in $SL(N, \mathbb{Z})$. By (2.61), the de Rham cohomology class of $\sum_{\mu \in \Lambda^*} \mu^* \delta_t$ is independent of $t$. The theorem follows from combining (2.62), (2.69) and (2.79). One can check that the arguments still go through if $E$ is not orientable. ∎

Recall from Section 1c that $n_{N,\mathbb{C}^N}^z$ denotes both an element of the (discrete) group cohomology $\mathrm{H}^{2N-1}(GL(N, \mathbb{C}); \mathbb{R})$ and the corresponding element of $\mathrm{H}^{2N-1}(BGL(N, \mathbb{C})_\delta; \mathbb{R})$.

**Theorem 2.22 :** *Let $i : GL(N, \mathbb{Z}) \to GL(N, \mathbb{C})$ be the natural inclusion. Then*

$$i^*\left(n_{N,\mathbb{C}^N}^z\right) = 0 \quad \text{in } \mathrm{H}^{2N-1}(GL(N, \mathbb{Z}); \mathbb{R}). \tag{2.81}$$

**Proof :** Put $\Gamma = GL(N, \mathbb{Z})$. It is known that there is a CW classifying space $B\Gamma$ with a finite number of cells in each dimension. Let $\mathcal{U}$ be the canonical flat $\mathbb{R}^N$-bundle over $B\Gamma$. For $K \gg 2N$, let $B\Gamma^K$ be the $K$-skeleton of $B\Gamma$. Let $B$ be a smooth connected compact manifold (possibly with boundary) which is homotopy-equivalent to $B\Gamma^K$. Let $\nu : B \to B\Gamma$ be the classifying map and put $E = \nu^*\mathcal{U}$. By Theorem 2.21, $n_N^z(\nabla^E) = 0$. Let $Bi : B\Gamma \to BGL(N, \mathbb{C})_\delta$ be the map on classifying spaces induced by $i$. By universality, $n_N^z(\nabla^E) = \nu^*(Bi)^*\left(n_{N,\mathbb{C}^N}^z\right)$. As $\nu$ is highly connected, it follows that $(Bi)^*\left(n_{N,\mathbb{C}^N}^z\right) = 0$, which is equivalent to the statement of the theorem. ∎

**Definition 2.23 :** If $N > 1$ and $s \in \mathbb{C}$, define $\phi(s) \in \Omega^{2N-2}(B)$ by

$$\phi(s) = -\int_0^\infty t^s \sum_{\mu \in \Lambda^*} \mu^* \epsilon_t \, dt. \tag{2.82}$$

Using (2.63) and (2.80), we see that $\phi(s)$ is well-defined and is a holomorphic function on $\mathbb{C}$.

**Theorem 2.24 :** *For* $\text{Re}(s) << 0$,

$$\phi(s) = 2^{-N-2s} \, \pi^{-\frac{N}{2}} \, \Gamma(N + \frac{1}{2} - s) \, \text{Vol}(E/\Lambda) \sum_{\substack{m \in \Lambda \\ m \neq 0}} (|m|)^{2s-2N-1} \, m^* \left( i_x \text{Vol} \cdot \int^B \widehat{\mathbf{x}} e^{-\frac{1}{8}\widehat{\omega}^2} \right). \quad (2.83)$$

**Proof :** This follows from (2.47) and (2.70). ∎

It follows that the right-hand-side of (2.83) must also have a holomorphic continuation to $\mathbb{C}$.

**Theorem 2.25 :** *If* $N > 1$ *then*

$$d\phi(0) = \frac{1}{2} \, \pi^{N-1} \, c_N^z \left( \nabla^E, h^E \right). \quad (2.84)$$

**Proof :** This follows from (2.61), (2.62), (2.79) and (2.82). ∎

## A - Results from [BLo] and vanishing of group cohomology classes

In this appendix we describe results from [BLo] on flat vector bundles and their direct images.

The appendix is organized as follows. In a) we establish some conventions. In b) we review results from [BLo] on finite-dimensional flat superconnections, their characteristic classes and their torsion forms. In c) we describe the relevant characteristic classes of flat vector bundles. In d) we review results from [BLo] on fiber bundles and the direct images of flat vector bundles. In e) we give a vanishing statement concerning the cohomology of a group which acts on a smooth compact manifold. We show how this implies Theorem 0.2.

### a) Conventions

We follow the notation of Section 1. For background information on superconnections we refer to [B], [BGV] and [Q]. Let $\varphi : \Omega(B) \to \Omega(B)$ be the linear map such that for all homogeneous $\omega \in \Omega(B)$,

$$(A.1) \qquad \varphi \, \omega = (2i\pi)^{-(\deg \omega)/2} \omega.$$

Let $E = E_+ \oplus E_-$ be a $\mathbb{Z}_2$-graded finite-dimensional vector bundle on $B$. Let $\tau$ be the involution of $E$ defining the $\mathbb{Z}_2$-grading, so that $\tau\big|_{E_\pm} = \pm I$. Then $\mathrm{End}(E)$ is a $\mathbb{Z}_2$-graded bundle of algebras over $B$, whose even (resp. odd) elements commute (resp. anticommute) with $\tau$. Given $a \in C^\infty(B; \mathrm{End}(E))$, we define its supertrace $\mathrm{Tr}_s[a] \in C^\infty(B)$ by

$$(A.2) \qquad \mathrm{Tr}_s[a] = \mathrm{Tr}[\tau a].$$

Given $\omega \in \Omega(B)$ and $a \in C^\infty(B; \mathrm{End}(E))$, put

$$(A.3) \qquad \mathrm{Tr}_s[\omega \cdot a] = \omega \, \mathrm{Tr}_s[a].$$

Then $\mathrm{Tr}_s$ extends to a linear map from $\Omega(B; \mathrm{End}(E))$ to $\Omega(B)$.

Given $\alpha, \alpha' \in \Omega(B; \mathrm{End}(E))$, we define their supercommutator $[\alpha, \alpha'] \in \Omega(B; \mathrm{End}(E))$ by

$$(A.4) \qquad [\alpha, \alpha'] = \alpha \alpha' - (-1)^{(\deg \alpha)(\deg \alpha')} \alpha' \alpha.$$

A basic fact is that $\mathrm{Tr}_s$ vanishes on supercommutators [Q].

Let $\nabla^E$ be a connection on $E$ which preserves the splitting $E = E_+ \oplus E_-$. Then $\nabla^E$ decomposes as $\nabla^E = \nabla^{E_+} \oplus \nabla^{E_-}$. Let $S$ be an odd element of $\Omega(B; \mathrm{End}(E))$. By definition, $\nabla^E + S$ gives a superconnection $A$ on $E$. That is, there is a $\mathbb{C}$-linear map

$$A : C^\infty(B; E) \to \Omega(B; E)$$

which is odd with respect to the total $\mathbb{Z}_2$-gradings and satisfies the Leibniz rule. We can extend $A$ to an odd $\mathbb{C}$-linear endomorphism of $\Omega(B; E)$. By definition, the curvature of $A$ is $A^2$, an even $C^\infty(B)$-linear endomorphism of $\Omega(B; E)$ which is given by multiplication by an even element of $\Omega(B; \operatorname{End}(E))$.

In what follows, we will say that a holomorphic function $f : \mathbb{C} \to \mathbb{C}$ is real if for all $a \in \mathbb{C}$, we have $f(\overline{a}) = \overline{f(a)}$.

**b) Characteristic classes and torsion forms of flat superconnections**

**Definition A.1 :** A superconnection $A'$ on $E$ is flat if its curvature vanishes, i.e. if $A'^2 = 0$.

Hereafter we assume that $A'$ is flat. Let $h^E$ be a Hermitian metric on $E$ such that $E_+$ and $E_-$ are orthogonal. Then there is a flat superconnection $A'^*$ on $E$ which is the adjoint of $A'$ with respect to $h^E$. Define an odd element of $\Omega(B; \operatorname{End}(E))$ by

$$(A.5) \qquad X = \frac{1}{2}(A'^* - A').$$

**Definition A.2 :** Let $f : \mathbb{C} \to \mathbb{C}$ be a holomorphic real odd function. Put

$$(A.6) \qquad f\left(A', h^E\right) = (2i\pi)^{1/2} \varphi \operatorname{Tr_s}[f(X)] \in \Omega(B).$$

**Theorem A.3 [BLo] :** *The differential form $f\left(A', h^E\right)$ is real, odd and closed. Its de Rham cohomology class is independent of $h^E$.*

We will denote the de Rham cohomology class of $f\left(A', h^E\right)$ by $f(A') \in \operatorname{H}^{\operatorname{odd}}(B; \mathbb{R})$.

Suppose now that $E = \bigoplus_{i=0}^{n} E^i$ is a $\mathbb{Z}$-graded complex vector bundle on $B$. Put

$$E_+ = \bigoplus_{i \text{ even}} E^i \quad, \quad E_- = \bigoplus_{i \text{ odd}} E^i.$$

Then $E = E_+ \oplus E_-$ is a $\mathbb{Z}_2$-graded vector bundle, to which we may apply the above formalism.

Let $A'$ be a superconnection on $E$. We can expand $A$ as

$$(A.7) \qquad A = \sum_{j \geq 0} A_j,$$

where $A_j$ is of partial degree $j$ with respect to the $\mathbb{Z}$-grading on $\Lambda\left(T^*B\right)$.

**Definition A.4 :** We say that $A'$ is of total degree 1 if

- $A'_1$ is a connection on $E$ which preserves the $\mathbb{Z}$-grading.
- For $j \in \mathbb{N} - \{1\}$, $A'_j$ is an element of $\Omega^j\left(B; \text{Hom}\left(E^\bullet, E^{\bullet+1-j}\right)\right)$.

In what follows, we will assume that $A'$ is a flat superconnection of total degree 1. Put

$$(A.8) \qquad v = A'_0 \ , \quad \nabla^E = A'_1.$$

Clearly $v \in C^\infty\left(B; \text{Hom}\left(E^\bullet, E^{\bullet+1}\right)\right)$. The flatness of $A'$ implies that

$$(A.9) \qquad v^2 = \left[\nabla^E, v\right] = \left(\nabla^E\right)^2 + [v, A'_2] = 0.$$

As $v^2 = 0$, we have a cochain complex of vector bundles

$$(A.10) \qquad (E, v) : 0 \longrightarrow E^0 \xrightarrow{v} E^1 \xrightarrow{v} \cdots \xrightarrow{v} E^n \longrightarrow 0.$$

**Definition A.5 :** For $b \in B$, let $H(E, v)_b = \bigoplus_{i=0}^{n} H^i(E, v)_b$ be the cohomology of the complex $(E, v)_b$.

Using (A.9), one can show that there a $\mathbb{Z}$-graded complex vector bundle $H(E, v)$ on $B$ whose fiber over $b \in B$ is $H(E, v)_b$, and a natural flat connection $\nabla^{H(E,v)}$ on $H(E, v)$.

Let $h^E$ be a Hermitian metric on $E$ such that the $E^i$'s are mutually orthogonal. Put $A'' = A'^*$, the adjoint superconnection to $A'$ with respect to $h^E$. Let $v^* \in C^\infty\left(B; \text{Hom}\left(E^\bullet, E^{\bullet-1}\right)\right)$ be the adjoint of $v$ with respect to $h^E$. From finite-dimensional Hodge theory, there is an isomorphism

$$(A.11) \qquad H(E, v) \cong \text{Ker}(v^*v + vv^*).$$

Being a subbundle of $E$, the vector bundle $\text{Ker}(v^*v + vv^*)$ inherits a Hermitian metric $h^{\text{Ker}}$ from the Hermitian metric $h^E$ on $E$. Let $h^{H(E,v)}$ denote the Hermitian metric on $H(E, v)$ obtained from $h^{\text{Ker}}$ via the isomorphism (A.11).

Let $N \in \text{End}(E)$ be the number operator of $E$, i.e. $N$ acts on $E^i$ by multiplication by $i$. Extend $N$ to an element of $C^\infty(B; \text{End}(E))$.

**Definition A.6 :** For $t > 0$, let $C'_t$ be the flat superconnection on $E$ of total degree 1 given by

$$(A.12) \qquad C'_t = t^{N/2} \, A' \, t^{-N/2}.$$

and let $C''_t$ be the flat superconnection on $E$ given by

$$(A.13) \qquad C''_t = t^{-N/2} \, A'' \, t^{N/2}.$$

The superconnections $C'_t$ and $C''_t$ are adjoint with respect to $h^E$. We have

$$(A.14) \qquad C'_t = \sum_{j \geq 0} t^{(1-j)/2} A'_j,$$

$$C''_t = \sum_{j \geq 0} t^{(1-j)/2} A''_j.$$

Define an odd element of $\Omega(B; \mathrm{End}(E))$ by

$$(A.15) \qquad D_t = \frac{1}{2} \left( C''_t - C'_t \right).$$

**Definition A.7 :** Define a real even differential form on $B$ by

$$(A.16) \qquad f^{\wedge}(C'_t, h^E) = \varphi \, \mathrm{Tr}_s \left[ \frac{N}{2} f'(D_t) \right] \in \Omega(B).$$

**Theorem A.8 [BLo] :** *One has*

$$(A.17) \qquad \frac{\partial}{\partial t} f(C'_t, h^E) = \frac{1}{t} df^{\wedge}(C'_t, h^E).$$

Let $d(H(E,v))$ be the constant integer-valued function on $B$

$$(A.18) \qquad d(H(E,v)) = \sum_{i=0}^{n} (-1)^i i \, \mathrm{rk}\left(H^i(E,v)\right).$$

In the rest of this section, we take $f(z) = z \, \exp(z^2)$.

**Theorem A.9 [BLo]** : $As\ t \to +\infty$,

(A.19)
$$f\left(C'_t, h^E\right) = f\left(\nabla^{H(E,v)}, h^{H(E,v)}\right) + \mathcal{O}\left(\frac{1}{\sqrt{t}}\right),$$
$$f^\wedge\left(C'_t, h^E\right) = d(H(E,v))\frac{f'(0)}{2} + \mathcal{O}\left(\frac{1}{\sqrt{t}}\right).$$

Now consider the special case when the vector bundle $E$ has not only a flat superconnection, but has a flat connection. Let

(A.20)
$$(E, v) : 0 \longrightarrow E^0 \xrightarrow{v} E^1 \xrightarrow{v} \cdots \xrightarrow{v} E^n \longrightarrow 0$$

be a flat complex of complex vector bundles. That is,

(A.21)
$$\nabla^E = \bigoplus_{i=0}^{n} \nabla^{E^i}$$

is a flat connection on $E = \bigoplus_{i=0}^{n} E^i$ and $v$ is a flat cochain map, meaning

(A.22)
$$\left(\nabla^E\right)^2 = 0\ ,\quad v^2 = 0\ ,\quad \left[\nabla^E, v\right] = 0.$$

For $t > 0$, put

(A.23)
$$C'_t = \sqrt{t}\ v + \nabla^E,$$
$$C''_t = \sqrt{t}\ v^* + \left(\nabla^E\right)^*.$$

Then $C'_t$ is a flat superconnection of total degree 1. Let $d(E)$ be the constant integer-valued function on $B$ given by

(A.24)
$$d(E) = \sum_{i=0}^{n} (-1)^i i \operatorname{rk}\left(E^i\right).$$

**Theorem A.10 [BLo]** : $As\ t \to 0$,

(A.25)
$$f\left(C'_t, h^E\right) = f\left(\nabla^E, h^E\right) + \mathcal{O}(t),$$
$$f^\wedge\left(C'_t, h^E\right) = d(E)\frac{f'(0)}{2} + \mathcal{O}(t).$$

**Theorem A.11 [BLo] :** *As elements of* $\mathrm{H}^{\mathrm{odd}}(B;\mathbb{R})$,

$$(A.26) \qquad f\left(\nabla^E\right) = f\left(\nabla^{H(E,v)}\right).$$

We now refine (A.26) to a statement about differential forms on $B$.

**Definition A.12 :** Define a real even differential form on $B$ by
$$(A.27)$$
$$T_f(A', h^E) = -\int_0^{+\infty} \left[ f^{\wedge}\left(C'_t, h^E\right) - d(H(E,v))\frac{f'(0)}{2} - [d(E) - d(H(E,v))]\frac{f'(\frac{i\sqrt{t}}{2})}{2} \right] \frac{dt}{t}.$$

**Remark A.13 :** By Theorems A.9 and A.10, the integrand in (A.27) is integrable. We will call $T_f(A', h^E)$ a torsion form.

**Theorem A.14 [BLo] :** *One has*

$$(A.28) \qquad dT_f(A', h^E) = f\left(\nabla^E, h^E\right) - f\left(\nabla^{H(E,v)}, h^{H(E,v)}\right).$$

**Proof :** This follows from Theorems A.8, A.9 and A.10. ∎

Upon passing to de Rham cohomology, Theorem A.14 implies Theorem A.11. Up to an overall multiplicative constant, the 0-form part of the torsion form $T_f(A', h^E)$ is the function which to a point $b \in B$ assigns the torsion of the cochain complex $(E, v)_b$, in the sense of [M, RS].

**c) Relationship to the characteristic classes of flat vector bundles**

We now relate the constructions of Section 1b and the preceding section. As in Section 1b, let $F$ be a complex vector bundle on $B$, endowed with a flat connection $\nabla^F$. We can consider $F$ to be a $\mathbb{Z}_2$-graded vector bundle with $F_+ = F$ and $F_- = 0$. Let $h^F$ be a Hermitian metric on $F$. In the rest of this section we will abbreviate $\omega\left(\nabla^F, h^F\right)$ by $\omega$.

Taking the flat superconnection $A'$ to be $\nabla^F$, the $X$ of (A.5) is given by

$$(A.29) \qquad X = \frac{\omega}{2}.$$

Let $f : \mathbb{C} \to \mathbb{C}$ be a holomorphic real odd function. Following Definition A.2, put

$$(A.30) \qquad f\left(\nabla^F, h^F\right) = (2i\pi)^{1/2} \varphi \operatorname{Tr}\left[f\left(\frac{\omega}{2}\right)\right] \in \Omega^{\mathrm{odd}}(B).$$

By Theorem A.3, $f\left(\nabla^F, h^F\right)$ is a real closed differential form on $B$ whose cohomology class $f\left(\nabla^F\right) \in \mathrm{H}^{\mathrm{odd}}(B;\mathbb{R})$ is independent of $h^F$.

If $j$ is a positive integer and $f(z) = z^{2j-1}$ then

$$(A.31) \qquad f(\nabla^F, h^F) = (2i\pi)^{-(j-1)}\, 2^{-(2j-1)}\, \mathrm{Tr}\left[\omega^{2j-1}\right].$$

In the notation of Section 1b, this equals $\frac{1}{j}\, n_j^z(\nabla^F, h^F)$.

### d) Fiber bundles

Let $Z \to M \xrightarrow{\pi} B$ be a smooth fiber bundle with connected base $B$ and connected closed fibers $Z_b = \pi^{-1}(b)$. Let $F$ be a flat complex vector bundle on $M$. Let $H\left(Z; F\big|_Z\right)$ denote the $\mathbb{Z}$-graded complex vector bundle on $B$ whose fiber over $b \in B$ is isomorphic to the cohomology group $\mathrm{H}^*(Z_b, F\big|_{Z_b})$. It has a canonical flat connection $\nabla^{H(Z;F|_Z)}$ which preserves the $\mathbb{Z}$-grading. Let $TZ$ be the vertical tangent bundle of the fiber bundle and let $o(TZ)$ be its orientation bundle, a flat real line bundle on $M$. Let $e(TZ) \in \mathrm{H}^{\dim(Z)}(M; o(TZ))$ be the Euler class of $TZ$. Put $f(z) = z\,\exp\left(z^2\right)$.

**Theorem A.15 [BLo]** : *One has an equality in* $\mathrm{H}^{\mathrm{odd}}(B;\mathbb{R})$:

$$(A.32) \qquad f\left(\nabla^{H(Z;F|_Z)}\right) = \int_Z e(TZ) \cdot f(\nabla^F).$$

In fact, one can refine (A.32) to a statement about differential forms on $B$. First, equip the fiber bundle with a horizontal distribution $T^H M$. Let $W$ be the infinite-dimensional $\mathbb{Z}$-graded vector bundle on $B$ whose fiber over $b \in B$ is isomorphic to $\Omega\left(Z_b; F\big|_{Z_b}\right)$. Then

$$(A.33) \qquad C^\infty(B; W) \simeq C^\infty(M; \Lambda\left(T^*Z\right) \otimes F)$$

and there is an isomorphism of $\mathbb{Z}$-graded vector spaces

$$(A.34) \qquad \Omega\left(M; F\right) \simeq \Omega\left(B; W\right).$$

Let $N$ be the number operator of $W$; it acts as multiplication by $i$ on $C^\infty(M; \Lambda^i(T^*Z) \otimes F)$.

The exterior differentiation operator $d^M$, acting on $\Omega\left(M; F\right)$, defines a flat superconnection on $W$ of total degree 1. In terms of the $\mathbb{Z}$-grading on $\Lambda\left(T^*B\right)$, $d^M$ can be decomposed as

$$(A.35) \qquad d^M = d^Z + \nabla^W + i_T,$$

where $d^Z$ is vertical exterior differentiation, $\nabla^W$ is a natural connection on $W$ which preserves the $\mathbb{Z}$-grading and $i_T$ is interior multiplication by the curvature $T$ of the fiber bundle, a $TZ$-valued horizontal 2-form on $M$. For $t > 0$, put

$$(A.36) \quad \begin{aligned} C'_t &= t^{N/2} \, d^M \, t^{-N/2} \\ &= \sqrt{t} \, d^Z + \nabla^W + \frac{1}{\sqrt{t}} \, i_T. \end{aligned}$$

Now equip the fiber bundle with a vertical Riemannian metric $g^{TZ}$ and the flat vector bundle $F$ with a Hermitian metric $h^F$. Then $W$ acquires an $L^2$-inner product $h^W$. There is a canonical metric-compatible connection $\nabla^{TZ}$ on $TZ$ [B, BGV]. The vector bundle $H\left(Z; F|_Z\right)$ acquires a Hermitian metric $h^{H(Z;F|_Z)}$ from Hodge theory. Let $C''_t$ be the adjoint superconnection to $C'_t$ with respect to $h^W$. That is,

$$(A.37) \quad C''_t = \sqrt{t} \, (d^Z)^* + (\nabla^W)^* - \frac{1}{\sqrt{t}} \, (T \wedge).$$

Define $D_t$, an odd element of $\Omega(B; \mathrm{End}(W))$, by

$$(A.38) \quad D_t = \frac{1}{2} \left( C''_t - C'_t \right).$$

For $t > 0$, define a real odd differential form on $B$ by

$$(A.39) \quad f\left(C'_t, h^W\right) = (2i\pi)^{1/2} \, \varphi \, \mathrm{Tr}_s \left[ f(D_t) \right]$$

and a real even differential form on $B$ by

$$(A.40) \quad f^\wedge \left(C'_t, h^W\right) = \varphi \, \mathrm{Tr}_s \left[ \frac{N}{2} f'(D_t) \right].$$

**Theorem A.16 [BLo] :** *For any $t > 0$,*

$$(A.41) \quad \frac{\partial}{\partial t} f\left(C'_t, h^W\right) = \frac{1}{t} df^\wedge \left(C'_t, h^W\right).$$

Put

$$(A.42) \quad \chi'(Z; F) = \sum_{i=0}^{\dim(Z)} (-1)^i i \, \mathrm{rk}\left( H^i(Z; F|_Z) \right),$$

an integer-valued constant function on $B$.

**Theorem A.17 [BLo] :** *As $t \to 0$,*

$$(A.43) \quad f\left(C'_t, h^W\right) = \int_Z e\left(TZ, \nabla^{TZ}\right) f\left(\nabla^F, h^F\right) + \mathcal{O}(t) \quad \text{if } \dim(Z) \text{ is even}$$

$$= \mathcal{O}\left(\sqrt{t}\right) \quad \text{if } \dim(Z) \text{ is odd,}$$

$$f^\wedge\left(C'_t, h^W\right) = \frac{1}{4}\dim(Z)\operatorname{rk}(F)\chi(Z) + \mathcal{O}(t) \quad \text{if } \dim(Z) \text{ is even,}$$

$$= \mathcal{O}\left(\sqrt{t}\right) \quad \text{if } \dim(Z) \text{ is odd.}$$

*As $t \to +\infty$,*

$$(A.44) \quad f\left(C'_t, h^W\right) = f\left(\nabla^{H(Z;F|_Z)}, h^{H(Z;F|_Z)}\right) + \mathcal{O}\left(\frac{1}{\sqrt{t}}\right),$$

$$f^\wedge\left(C'_t, h^W\right) = \frac{\chi'(Z;F)}{2} + \mathcal{O}\left(\frac{1}{\sqrt{t}}\right).$$

**Definition A.18 :** The analytic torsion form $\mathcal{T}\left(T^H M, g^{TZ}, h^F\right)$, a real even differential form on $B$, is given by

$$(A.45) \quad \mathcal{T}\left(T^H M, g^{TZ}, h^F\right) = -\int_0^{+\infty} \left[f^\wedge\left(C'_t, h^W\right) - \frac{\chi'(Z;F)}{2}f'(0)\right.$$

$$\left. - \left(\frac{\dim(Z)\operatorname{rk}(F)\chi(Z)}{4} - \frac{\chi'(Z;F)}{2}\right)f'(\frac{i\sqrt{t}}{2})\right]\frac{dt}{t}.$$

**Remark A.19 :** It follows from Theorem A.17 that the integrand of (A.45) is integrable.

**Theorem A.20 [BLo] :** *One has*

$$(A.46) \quad d\mathcal{T}\left(T^H M, g^{TZ}, h^F\right) = \int_Z e\left(TZ, \nabla^{TZ}\right) f\left(\nabla^F, h^F\right) - f\left(\nabla^{H(Z;F|_Z)}, h^{H(Z;F|_Z)}\right).$$

**Proof :** This follows from Theorems A.16 and A.17. ∎

Upon passing to de Rham cohomology, Theorem A.20 implies Theorem A.15.

**Remark A.21 :** One can extend the results of this section to the case when the fiber has boundary, by using the doubling trick of [LR, Section IX].

### e) Vanishing of group cohomology classes

We use the notation of Section 1c. Let $\Gamma$ be a discrete group. Suppose that there is a $CW$ complex which is a $K(\Gamma,1)$ space, with a finite number of cells in each dimension. We will denote this $CW$-complex by $B\Gamma$. Let $Z$ be a connected closed smooth manifold and let $\Gamma$ act on $Z$ by a homomorphism $\rho : \Gamma \to \mathrm{Diff}(Z)$. For each integer $p \in [0, \dim(Z)]$, there is an induced representation $r_p : \Gamma \to GL(\mathrm{H}^p(Z;\mathbb{C}))$. If $j$ is a positive integer, we can pullback the group cohomology class $n^z_{j,\mathrm{H}^p(Z;\mathbb{C})}$ under $r_p$ to obtain $r_p^*\left(n^z_{j,\mathrm{H}^p(Z;\mathbb{C})}\right) \in \mathrm{H}^{2j-1}(\Gamma;\mathbb{R})$.

**Theorem A.22 :** *For any positive integer $j$, we have an equality in $\mathrm{H}^{2j-1}(\Gamma;\mathbb{R})$:*

$$(A.47) \qquad \sum_{p=0}^{\dim(Z)} (-1)^p \, r_p^*\left(n^z_{j,\mathrm{H}^p(Z;\mathbb{C})}\right) = 0.$$

**Proof :** For $K \gg 2j-1$, let $B\Gamma^K$ be the $K$-skeleton of $B\Gamma$. Let $B$ be a smooth connected compact manifold (possibly with boundary) which is homotopy equivalent to $B\Gamma^K$. Let $\widetilde{B}$ be the universal cover of $B$. Put $M = \widetilde{B} \times_\rho Z$, the total space of a fiber bundle over $B$ with fiber $Z$. Let $F$ be the trivial flat complex line bundle on $M$. Then $f(\nabla^F) = 0$ and Theorem A.15 implies that

$$(A.48) \qquad f\left(\nabla^{H(Z;F|_Z)}\right) = 0.$$

The term of degree $2j-1$ in (A.48) gives

$$(A.49) \qquad n^z_j\left(\nabla^{H(Z;F|_Z)}\right) = 0.$$

By construction, the $\mathbb{Z}$-graded vector bundle $H(Z;F|_Z)$ on $B$ is such that

$$(A.50) \qquad H^p(Z;F|_Z) = \widetilde{B} \times_{r_p} \mathrm{H}^p(Z;\mathbb{C}).$$

Let $h : B \to B\Gamma$ be the classifying map and let $\nu_p$ be the composite

$$(A.51) \qquad B \xrightarrow{h} B\Gamma \xrightarrow{Br_p} BGL(\mathrm{H}^p(Z;\mathbb{C}))_\delta.$$

Using the discussion of Section 1c on the $n^z_j$-classes and (A.50), it follows that

$$(A.52) \qquad \sum_{p=0}^{\dim(Z)} (-1)^p \, \nu_p^*\left(n^z_{j,\mathrm{H}^p(Z;\mathbb{C})}\right) = 0.$$

With our assumptions on $K$ and $B$, the map $h$ induces an isomorphism $h^* : \mathrm{H}^{2j-1}(B\Gamma; \mathbb{R}) \to \mathrm{H}^{2j-1}(B; \mathbb{R})$. Thus we have an equality in $\mathrm{H}^{2j-1}(B\Gamma; \mathbb{R})$:

$$(A.53) \qquad \sum_{p=0}^{\dim(Z)} (-1)^p \, (Br_p)^* \left( n^z_{j, \mathrm{H}^p(Z;\mathbb{C})} \right) = 0.$$

Using the isomorphism between $\mathrm{H}^{2j-1}(B\Gamma; \mathbb{R})$ and $\mathrm{H}^{2j-1}(\Gamma; \mathbb{R})$, the theorem follows. ∎

**Remark A.23 :** One can give a more direct proof of Theorem A.22 in the case $j = 1$. If $j = 1$, equation (A.47) says that for all $\gamma \in \Gamma$,

$$(A.54) \qquad \sum_{p=0}^{\dim(Z)} (-1)^p \, \ln \left| \det r_p(\gamma) \big|_{\mathrm{H}^p(Z;\mathbb{C})} \right| = 0.$$

Let $T^p$ denote the torsion subgroup of $\mathrm{H}^p(Z; \mathbb{Z})$. Then $\mathrm{H}^p(Z; \mathbb{C}) \cong (\mathrm{H}^p(Z; \mathbb{Z})/T^p) \otimes \mathbb{C}$. As $\rho(\gamma)$ is a diffeomorphism of $Z$, it acts as an automorphism of the lattice $\mathrm{H}^p(Z; \mathbb{Z})/T^p \subset \mathrm{H}^p(Z; \mathbb{C})$. Thus $\det r_p(\gamma)\big|_{\mathrm{H}^p(Z;\mathbb{C})} = \pm 1$. Equation (A.54) follows.

We now apply Theorem A.22 to the special case when $N$ is a positive odd integer, $Z = T^N$ and $\Gamma = GL(N, \mathbb{Z})$, acting linearly on the torus.

**Theorem A.24 :** *Let $i : GL(N, \mathbb{Z}) \to GL(N, \mathbb{C})$ be the natural inclusion. Then*

$$(A.55) \qquad i^* \left( n^z_{N, \mathbb{C}^N} \right) = 0 \quad \text{in } \mathrm{H}^{2N-1}(GL(N, \mathbb{Z}); \mathbb{R}).$$

**Proof :** Put $\Gamma = GL(N, \mathbb{Z})$. It is known that there is a CW classifying space $B\Gamma$ with a finite number of cells in each dimension. Writing $T^N$ as $\mathbb{R}^N/\mathbb{Z}^N$, we have $(\mathbb{R}^N)^* \otimes \mathbb{C} \cong \overline{\mathbb{C}^N}^*$ and $\mathrm{H}^*(T^N; \mathbb{C}) \cong \Lambda\left(\overline{\mathbb{C}^N}^*\right)$. From Corollary 1.12,

$$(A.56) \qquad \mathrm{ch}^z_{\Lambda(\overline{\mathbb{C}^N}^*)} = \frac{1}{N} \, n^z_{N, \mathbb{C}^N}.$$

From the definition of $\mathrm{ch}^z$, the term of degree $2N - 1$ of $\mathrm{ch}^z_{\Lambda(\overline{\mathbb{C}^N}^*)}$ is proportionate to $n^z_{N, \Lambda(\overline{\mathbb{C}^N}^*)}$. The theorem now follows from applying (A.47) in the case $j = N$. ∎

## B - Fourier analysis on torus bundles

In this appendix we describe how in the special case of a flat torus bundle, the results of Appendix A become equivalent to the results of Section 2.

The appendix is organized as follows. In a) we review the relationship between supertraces and the Berezin integral. Given a flat rank-$N$ vector bundle $\mathcal{V}$, in b) we define its leafwise topology $\mathcal{V}_\mathcal{F}$ and use the formalism of Appendix A to define a closed form $\widetilde{\delta}_t \in \Omega^{2N-1}(\mathcal{V}_\mathcal{F})$. We also construct the transgressing form $\widetilde{\epsilon}_t \in \Omega^{2N-1}(\mathcal{V}_\mathcal{F})$, and give the relationship between $(\widetilde{\delta}_t, \widetilde{\epsilon}_t)$ and the forms $(\delta_t, \epsilon_t)$ of Section 2. In c) we use Fourier analysis on the fibers of a torus bundle to show how the results of Appendix A become equivalent to those of Section 2.

### a) Supertraces and the Berezin integral

Let $V$ be a Hermitian vector space of dimension $N$ with orthonormal basis $\{e_k\}_{k=1}^N$. Given $X \in V$, define operators on $\Lambda(V)$ by

$$(B.1) \qquad c(X) = (X \wedge) - i(X),$$
$$\widehat{c}(X) = (X \wedge) + i(X).$$

Then for $X, Y \in V$,

$$(B.2) \qquad c(X)c(Y) + c(Y)c(X) = -2 \langle X, Y \rangle,$$
$$\widehat{c}(X)\widehat{c}(Y) + \widehat{c}(Y)\widehat{c}(X) = 2 \langle X, Y \rangle,$$
$$c(X)\widehat{c}(Y) + \widehat{c}(Y)c(X) = 0.$$

Thus $c$ and $\widehat{c}$ generate two graded-commuting Clifford algebras. Put $c_i = c(e_i)$ and $\widehat{c}_i = \widehat{c}(e_i)$. Among the monomials in the $c_i$'s and $\widehat{c}_i$'s with less than or equal to $N$ factors of each, the only nonzero supertrace occurs as

$$(B.3) \qquad \mathrm{Tr}_s [c_1 \ldots c_N \widehat{c}_1 \ldots \widehat{c}_N] = (-1)^{\frac{N(N+1)}{2}} 2^N.$$

We now relate the above supertrace to the Berezin integral of (2.6). Let $\widehat{V}$ be another copy of $V$. Let $M \in \mathrm{End}(V \oplus \widehat{V})$ be a skew-symmetric endomorphism. Let $\det^{1/2}\left(\frac{\sin(M)}{M}\right) \in \mathbb{C}$ be the square-root of $\det\left(\frac{\sin(M)}{M}\right)$ which extends to a holomorphic function on the space of skew-symmetric endomorphisms, with value 1 at $M = 0$. As notation, we write $C$ and $\Psi$ for the $(2n)$-vectors $\begin{pmatrix} c \\ \widehat{c} \end{pmatrix}$ and $\begin{pmatrix} \psi \\ \widehat{\psi} \end{pmatrix}$, respectively. Let $A$ be a graded-commutative superalgebra. Let $J = \begin{pmatrix} j \\ \widehat{j} \end{pmatrix}$ be a $(2N)$-vector of odd elements of $A$.

The next theorem is a consequence of [MQ, (2.13)].

**Theorem B.1 :** *One has an identity in $A$ :*

$$(B.4) \quad \text{Tr}_s \left[ e^{\frac{1}{2} \langle C, MC \rangle + \langle J, C \rangle} \right] = (-1)^{\frac{N(N+1)}{2}} 2^N \det^{1/2} \left( \frac{\sin(M)}{M} \right) \int^B e^{\frac{1}{2} \langle \Psi, M\Psi \rangle + \langle J, \Psi \rangle}.$$

One can extend (B.4) to allow $M$ to have entries which are even elements of $A$, provided that the entries are nilpotent or that $A$ is a Banach algebra.

**b) Thom-like forms II**

Let $\mathcal{V}$ be a complex rank-$N$ vector bundle over a smooth connected manifold $B$. Suppose that $\mathcal{V}$ has a flat connection $\nabla^{\mathcal{V}}$. Let $\nabla^{\Lambda(\mathcal{V})}$ be the induced flat connection on the $\mathbb{Z}$-graded vector bundle $\Lambda(\mathcal{V})$. The tangent vectors of $\mathcal{V}$ which are horizontal for $\nabla^{\mathcal{V}}$ define an integrable distribution on $\mathcal{V}$, and hence a foliation $\mathcal{F}$ which is transverse to the fibers of $\mathcal{V}$. Let $\mathcal{V}_{\mathcal{F}}$ denote the total space of $\mathcal{V}$ with the leaf topology, a basis of which is given by the connected components of intersections $S \cap L$ of of open sets $S$ in $V$ with leaves $L$ of $\mathcal{F}$ [L, p. 2]. The connected components of $\mathcal{V}_{\mathcal{F}}$ are exactly the leaves of $\mathcal{F}$. (For example, if $B$ is a point then $\mathcal{V}_{\mathcal{F}}$ is $\mathbb{C}^N$ with the discrete topology.) In particular, if $U$ is a contractible open set in $B$ and $\lambda : U \to \mathcal{V}_{\mathcal{F}}$ is a continuous section of $\mathcal{V}_{\mathcal{F}}$ over $U$ then $\lambda$ is automatically a flat section. Note that $\mathcal{V}_{\mathcal{F}}$ is a non-second-countable manifold whose dimension is that of $B$.

Let $\pi : \mathcal{V}_{\mathcal{F}} \to B$ be the projection map. Put $\mathcal{E} = \Lambda(\pi^*\mathcal{V})$, a $\mathbb{Z}$-graded vector bundle on $\mathcal{V}_{\mathcal{F}}$. It is equipped with the flat connection $\nabla^{\mathcal{E}} = \pi^* \nabla^{\Lambda(\mathcal{V})}$.

**Definition B.2 :** The operator $(\mathbf{x} \wedge) \in C^{\infty}(\mathcal{V}_{\mathcal{F}}; \text{Hom}(\mathcal{E}^{\bullet}, \mathcal{E}^{\bullet+1}))$ acts at $x \in \mathcal{V}_{\mathcal{F}}$ as exterior multiplication by $x$ on $\mathcal{E}_x = \Lambda(\pi^*\mathcal{V})_x = \Lambda(\mathcal{V}_{\pi(x)})$.

**Definition B.3 :** The superconnection $A'$ on $\mathcal{E}$ is given by

$$(B.5) \quad A' = \sqrt{-1} \, (\mathbf{x} \wedge) + \nabla^{\mathcal{E}}.$$

Clearly $A'$ has total degree 1.

**Theorem B.4 :** *The superconnection $A'$ is flat.*

**Proof :** Clearly $(\mathbf{x} \wedge)^2 = (\nabla^{\mathcal{E}})^2 = 0$. Given $b \in B$, let $U$ be a contractible neighborhood of $b$. We can trivialize the vector bundle $\mathcal{V}$ over $U$ so that if $\sigma : U \to \mathbb{C}^N$ is a section of $\mathcal{V}|_U$ then $\nabla^\mathcal{V}(\sigma) = d\sigma$. Now $\pi^{-1}(U) \cong U \times \mathbb{C}^N$, and the topology on $U \times \mathbb{C}^N$ which is induced from $\mathcal{V}_\mathcal{F}$ is the product of the topology of $U$ with the discrete topology on $\mathbb{C}^N$. The operator $d$ on $\Omega\left(\pi^{-1}(U)\right)$ is effectively exterior differentiation in the $U$-direction on $U \times \mathbb{C}^N$, and we can write it as $d_U$. We have

$$(B.6) \qquad \Lambda\left(\pi^*\mathcal{V}\right)\big|_{\pi^{-1}(U)} \cong \left(U \times \mathbb{C}^N\right) \times \Lambda\left(\mathbb{C}^N\right).$$

Writing a local section $s : U \times \mathbb{C}^N \to \Lambda\left(\mathbb{C}^N\right)$ of $\mathcal{E}$ as $s(u, x)$, we have

$$(B.7) \qquad \begin{aligned}(\mathbf{x} \wedge s)(u, x) &= x \wedge s(u, x), \\ (\nabla^\mathcal{E} s)(u, x) &= d_U s(u, x).\end{aligned}$$

It is now clear that

$$(B.8) \qquad \left[\nabla^\mathcal{E}, (\mathbf{x} \wedge)\right] = 0.$$

The theorem follows. ∎

**Remark B.5 :** Theorem B.4 would not be true if we used the ordinary topology on $\mathcal{V}$. This can be seen in the case when $B$ is a point.

Let $h^\mathcal{V}$ be a Hermitian metric on $\mathcal{V}$. There is an induced Hermitian metric $h^\mathcal{E}$ on $\mathcal{E}$.

**Definition B.6 :** Let $i_\mathbf{x} \in C^\infty(\mathcal{V}_\mathcal{F}; \text{Hom}(\mathcal{E}^\bullet, \mathcal{E}^{\bullet-1}))$ be the operator which acts at $x \in \mathcal{V}_\mathcal{F}$ as interior multiplication by $x$ on $\mathcal{E}_x = \Lambda\left(\pi^*\mathcal{V}\right)_x = \Lambda\left(\mathcal{V}_{\pi(x)}\right)$.

Define $\omega\left(\nabla^\mathcal{V}, h^\mathcal{V}\right) \in \Omega^1(B; \text{End}(\mathcal{V}))$ as in (1.5). In the rest of this section we will abbreviate $\omega\left(\nabla^\mathcal{V}, h^\mathcal{V}\right)$ by $\omega$. We can lift $\omega$ to an operator $\pi^*\omega \in \Omega^1(\mathcal{V}_\mathcal{F}; \text{End}(\pi^*\mathcal{V}))$, and extend it to an operator $\Lambda\left(\pi^*\omega\right) \in \Omega^1(\mathcal{V}_\mathcal{F}; \text{End}(\Lambda(\pi^*\mathcal{V})))$.

Let $A'^*$ be the adjoint superconnection to $A'$ with respect to $h^\mathcal{E}$.

**Theorem B.7 :** *One has*

$$(B.9) \qquad A'^* = -\sqrt{-1}\, i_\mathbf{x} + \nabla^\mathcal{E} + \Lambda\left(\pi^*\omega\right).$$

**Proof :** Clearly $(\mathbf{x} \wedge)^* = i_\mathbf{x}$. In terms of the local trivializations of the proof of Theorem B.4, on $U$ we have

$$(B.10) \qquad \omega = \left(h^\mathcal{V}\right)^{-1}\left(dh^\mathcal{V}\right).$$

On $U \times \mathbb{C}^N$, we have

(B.11) $$\left(\nabla^{\mathcal{E}}\right)^* = d_U + \left(h^{\mathcal{E}}\right)^{-1} \left(d_U h^{\mathcal{E}}\right) = \nabla^{\mathcal{E}} + \Lambda\left(\pi^*\omega\right).$$

The theorem follows. ∎

Following the formalism of Appendix A, for $t > 0$ put

(B.12) $$C'_t = \sqrt{-1}\, \sqrt{t}\, (\mathbf{x} \wedge) + \nabla^{\mathcal{E}},$$
$$C''_t = -\sqrt{-1}\, \sqrt{t}\, i_{\mathbf{x}} + \nabla^{\mathcal{E}} + \Lambda\left(\pi^*\omega\right).$$

Then the $D_t$ of (A.15) is given by

(B.13) $$D_t = -\frac{1}{2}\,\sqrt{-1}\,\sqrt{t}\,\widehat{c}(\mathbf{x}) + \frac{1}{2}\Lambda\left(\pi^*\omega\right).$$

**Definition B.8 :** For $t > 0$, define $\widetilde{\delta}_t \in \Omega^{\text{odd}}(\mathcal{V}_{\mathcal{F}})$ and $\widetilde{\epsilon}_t \in \Omega^{\text{even}}(\mathcal{V}_{\mathcal{F}})$ by

(B.14) $$\widetilde{\delta}_t = f\left(C'_t, h^{\mathcal{E}}\right), \quad \widetilde{\epsilon}_t = \frac{1}{t}\, f^{\wedge}\left(C'_t, h^{\mathcal{E}}\right).$$

As a consequence of Theorems A.3 and A.8, we have that $\widetilde{\delta}_t$ is closed and

(B.15) $$\frac{\partial \widetilde{\delta}_t}{\partial t} = d\widetilde{\epsilon}_t.$$

**Definition B.9 :** Define $\psi \in \Omega^{\text{even}}(\mathcal{V}_{\mathcal{F}})$ by the right-hand-side of (A.27).

Let $s : B \to \mathcal{V}_{\mathcal{F}}$ be the zero-section.

**Theorem B.10 :** *On $\mathcal{V}_{\mathcal{F}} \setminus \text{Im}(s)$, the differential form $f\left(\nabla^{\mathcal{E}}, h^{\mathcal{E}}\right)$ is exact and*

(B.16) $$f\left(\nabla^{\mathcal{E}}, h^{\mathcal{E}}\right) = d\psi.$$

*In addition,*

(B.17) $$s^* f\left(\nabla^{\mathcal{E}}, h^{\mathcal{E}}\right) = f\left(\nabla^{\Lambda(\mathcal{V})}, h^{\Lambda(\mathcal{V})}\right).$$

**Proof :** In our case, the complex (A.20) becomes

$$(B.18) \qquad (\mathcal{E}, \mathbf{x} \wedge) : 0 \longrightarrow \Lambda^0(\pi^*\mathcal{V}) \xrightarrow{\mathbf{x}\wedge} \Lambda^1(\pi^*\mathcal{V}) \xrightarrow{\mathbf{x}\wedge} \cdots \xrightarrow{\mathbf{x}\wedge} \Lambda^n(\pi^*\mathcal{V}) \longrightarrow 0.$$

It is acyclic on $\mathcal{V}_\mathcal{F} \setminus \mathrm{Im}(s)$, and so (B.16) follows from Theorem A.14. Equation (B.17) follows from the naturality of the constructions. ∎

We now relate $\widetilde{\delta}_t$ and $\widetilde{\epsilon}_t$ to the forms $\delta_t$ and $\epsilon_t$ of Definition 2.8.

**Theorem B.11 :** *Let $I : \mathcal{V}_\mathcal{F} \to \mathcal{V}$ be the identity map. Then*

$$(B.19) \qquad \widetilde{\delta}_t = (-1)^{\frac{N-1}{2}} \, 2 \, \pi^{-(N-1)} \, I^* \delta_{\frac{t}{4}}$$

*and*

$$(B.20) \qquad \widetilde{\epsilon}_t = (-1)^{\frac{N-1}{2}} \, \frac{1}{2} \, \pi^{-(N-1)} \, I^* \epsilon_{\frac{t}{4}}.$$

**Proof :** Let $z$ be an auxiliary odd variable, satisfying $z^2 = 0$. Define $\mathrm{Tr}_z$ as in (1.13). From (A.6), we have

$$(B.21) \qquad f\left(C'_t, h^\mathcal{E}\right) = (2i\pi)^{1/2} \varphi \, \mathrm{Tr}_\mathrm{s} \left[ D_t e^{D_t^2} \right]$$
$$= (2i\pi)^{1/2} \varphi \, \mathrm{Tr}_\mathrm{z} \left[ \mathrm{Tr}_\mathrm{s} \left[ e^{D_t^2 + zD_t} \right] \right].$$

As $\omega$ is a symmetric one-form, we have

$$(B.22) \qquad \Lambda(\omega) = \sum_{kl} \omega_{kl} \, (e_k \wedge) \, i(e_l) = \frac{1}{4} \sum_{kl} \omega_{kl} \, (\widehat{c}_k + c_k)(\widehat{c}_l - c_l)$$
$$= -\frac{1}{2} \sum_{k,l} \omega_{kl} \widehat{c}_k c_l = \frac{1}{2} \sum_{k,l} \widehat{c}_k \omega_{kl} c_l = \frac{1}{2} \langle \widehat{c}, \omega c \rangle.$$

For simplicity, in the rest of the proof we write $\omega$ in place of $\pi^*\omega$. Then we can write the $D_t$ of (B.13) as

$$(B.23) \qquad D_t = -\frac{i\sqrt{t}}{2} \left\langle \widehat{c}, \mathbf{x} + \frac{i}{2\sqrt{t}} \, \omega c \right\rangle.$$

One can check that

$$(B.24) \qquad D_t^2 + zD_t = -\frac{t}{4} \left| \mathbf{x} + \frac{i}{2\sqrt{t}} \, \omega c + \frac{i}{\sqrt{t}} \, z\widehat{c} \right|^2 + \frac{1}{16} \langle \widehat{c}, \omega^2 \widehat{c} \rangle.$$

From (B.4), we have

(B.25)
$$\operatorname{Tr}_z\left[\operatorname{Tr}_s\left[e^{D_t^2+zD_t}\right]\right] =$$
$$(-1)^{\frac{N(N+1)}{2}} 2^N \operatorname{Tr}_z\left[\det{}^{1/2}\left(\frac{\sin(M)}{M}\right)\int^B e^{-\frac{t}{4}\left|\mathbf{x}+\frac{i}{2\sqrt{t}}\omega\psi+\frac{i}{\sqrt{t}}z\widehat{\psi}\right|^2+\frac{1}{16}\left\langle\widehat{\psi},\omega^2\widehat{\psi}\right\rangle}\right],$$

where the matrix $M$ is given by

(B.26)
$$M = \begin{pmatrix} -\frac{\omega^2}{8} & \frac{z\omega}{4} \\ -\frac{z\omega}{4} & \frac{\omega^2}{8} \end{pmatrix}.$$

Due to the nature of the operator $\operatorname{Tr}_z$, we can consider it to be acting on either the first or second factor of the term in brackets in the right-hand-side of (B.25). As the $z$ in $M$ occurs in the off-diagonal entries,

(B.27)
$$\operatorname{Tr}_z\left[\det{}^{1/2}\left(\frac{\sin(M)}{M}\right)\right] = 0.$$

Thus

(B.28)
$$\operatorname{Tr}_z\left[\operatorname{Tr}_s\left[e^{D_t^2+zD_t}\right]\right] =$$
$$(-1)^{\frac{N(N+1)}{2}} 2^N \det{}^{1/2}\left(\frac{\sin(M_0)}{M_0}\right)\operatorname{Tr}_z\left[\int^B e^{-\frac{t}{4}\left|\mathbf{x}+\frac{i}{2\sqrt{t}}\omega\psi+\frac{i}{\sqrt{t}}z\widehat{\psi}\right|^2+\frac{1}{16}\left\langle\widehat{\psi},\omega^2\widehat{\psi}\right\rangle}\right],$$

where the matrix $M_0$ is given by

(B.29)
$$M_0 = \begin{pmatrix} -\frac{\omega^2}{8} & 0 \\ 0 & \frac{\omega^2}{8} \end{pmatrix}.$$

By Lemma 1.3, $\det^{1/2}\left(\frac{\sin(M_0)}{M_0}\right) = 1$. Thus

(B.30) $\operatorname{Tr}_z\left[\operatorname{Tr}_s\left[e^{D_t^2+zD_t}\right]\right]$
$$= (-1)^{\frac{N(N+1)}{2}} 2^N \operatorname{Tr}_z\left[\int^B e^{-\frac{t}{4}\left|\mathbf{x}+\frac{i}{2\sqrt{t}}\omega\psi+\frac{i}{\sqrt{t}}z\widehat{\psi}\right|^2+\frac{1}{16}\left\langle\widehat{\psi},\omega^2\widehat{\psi}\right\rangle}\right]$$
$$= (-1)^{\frac{N(N-1)}{2}} 2^N \operatorname{Tr}_z\left[\int^B e^{-\frac{t}{4}\left|\mathbf{x}+\frac{1}{2\sqrt{t}}\omega\psi+\frac{1}{\sqrt{t}}z\widehat{\psi}\right|^2-\frac{1}{16}\left\langle\widehat{\psi},\omega^2\widehat{\psi}\right\rangle}\right].$$

On the other hand, from (2.73),

$$(B.31) \qquad I^*\delta_{\frac{t}{4}} = \text{Tr}_z\left[\int^B e^{-\frac{t}{4}\left|\mathbf{x}+\frac{1}{2\sqrt{t}}\omega\psi+\frac{1}{\sqrt{t}}z\widehat{\psi}\right|^2 - \frac{1}{16}\left\langle\widehat{\psi},\omega^2\widehat{\psi}\right\rangle}\right].$$

Equation (B.19) follows from combining (B.21), (B.30) and (B.31). Equation (B.20) follows similarly. ∎

**Remark B.12 :** It would have been more natural to phrase the results of Section 2b in terms of the forms $I^*\delta_t$ and $I^*\epsilon_t$. This is because it is only the flat structure of $\mathcal{V}$ which counts, not its topological structure.

### c) Fourier decomposition

We follow the geometric setup of Section 2d. Put $M = E/\Lambda$. Then $M$ is the total space of a fiber bundle over $B$ with fiber $Z = T^N$. There is a horizontal distribution $T^H M$ on $M$ given by pushing forward the horizontal vectors for $\nabla^E$ under the quotient map $E \to M$. Let $F$ be the trivial complex line bundle on $M$, with the standard Hermitian metric $h^F$. Let $h^E$ be an inner-product on $E$. We assume that the induced volume forms on the fibers are preserved by $\nabla^E$. Hereafter, we write $\omega_{E^*}$ for $\omega\left(\nabla^{E^*}, h^{E^*}\right)$. The metric $h^E$ also induces a Riemannian metric $g^{TZ}$ on the vertical tangent bundle $TZ$.

Let $\overline{W}$ be the $\mathbb{Z}$-graded Hilbert bundle on $B$ whose fiber over $b \in B$ is isomorphic to $L^2(Z_b; \Lambda(T^*Z_b))$. As the torus $Z_b$ has trivial tangent bundle, we have isomorphisms

$$(B.32) \qquad C^\infty(Z_b; \Lambda(T^*Z_b)) \cong \Lambda(E_b^*) \otimes C^\infty(Z_b),$$
$$L^2(Z_b; \Lambda(T^*Z_b)) \cong \Lambda(E_b^*) \otimes L^2(Z_b).$$

Given $\mu_b$ in the lattice $\Lambda_b^*$, there is a well-defined function $e^{\sqrt{-1}\langle\mu_b,\cdot\rangle} \in C^\infty(Z_b)$. By Fourier analysis, for each $b \in B$ there is an orthogonal decomposition

$$(B.33) \qquad \overline{W}_b = \bigoplus_{\mu_b \in \Lambda_b^*} \Lambda(E_b^*) \otimes \mathbb{C} e^{\sqrt{-1}\langle\mu_b,\cdot\rangle}.$$

If $U$ is a contractible open subset of $B$ then the orthogonal decompositions of $\{\overline{W}_b\}_{b\in U}$ piece together to give an orthogonal decomposition

$$(B.34) \qquad L^2(U;\overline{W}|_U) = \bigoplus_{\mu \in C^\infty(U;\Lambda^*)} L^2(U;\Lambda(E^*|_U)) \otimes \mathbb{C} e^{\sqrt{-1}\langle\mu,\cdot\rangle}.$$

**Theorem B.13 :** *With respect to the orthogonal decomposition* (B.34), *the superconnection* $d^M$ *splits as*

$$(B.35) \qquad d^M = \bigoplus_{\mu \in C^\infty(U; \Lambda^*)} \left( \sqrt{-1} \, (\mu \wedge) + \nabla^{\Lambda(E^*|_U)} \right) \otimes I.$$

**Proof :** This follows from (A.35), but we will give a direct proof. We can trivialize $E$ over $U$ as $E = U \times \mathbb{C}^N$ so that if $\sigma : U \to \mathbb{C}^N$ is a section of $E|_U$ then $\nabla^E(\sigma) = d\sigma$. With respect to this trivialization, $M|_U = U \times T^N$. There is an isomorphism

$$(B.36) \qquad \alpha : \Omega(U) \, \widehat{\otimes} \, \Lambda\left((\mathbb{C}^N)^*\right) \otimes C^\infty(T^N) \to \Omega(U) \, \widehat{\otimes} \, \Omega\left(T^N\right).$$

Acting on $\Omega(U) \, \widehat{\otimes} \, \Omega\left(T^N\right)$, we have $d^M = (I \, \widehat{\otimes} \, d_{T^N}) + (d_U \, \widehat{\otimes} \, I)$. Then in terms of our trivializations, the superconnection $d^M$ can be written as $\alpha^{-1}\left((I \, \widehat{\otimes} \, d_{T^N}) + (d_U \, \widehat{\otimes} \, I)\right) \alpha$.

Given $f \in C^\infty(U)$, $s \in \Lambda\left((\mathbb{C}^N)^*\right)$ and $\mu \in C^\infty(U; \Lambda^*)$, we have

$$(B.37) \qquad \alpha \left( f \otimes s \otimes e^{\sqrt{-1} \langle \mu, \cdot \rangle} \right) = f \otimes s e^{\sqrt{-1} \langle \mu, \cdot \rangle},$$

$$(B.38) \qquad \left(I \, \widehat{\otimes} \, d_{T^N}\right) \left( f \otimes s e^{\sqrt{-1} \langle \mu, \cdot \rangle} \right) = f \otimes \left( \sqrt{-1} \, \mu \wedge s \, e^{\sqrt{-1} \langle \mu, \cdot \rangle} \right)$$

and

$$(B.39) \qquad \left(d_U \, \widehat{\otimes} \, I\right) \left( f \otimes s e^{\sqrt{-1} \langle \mu, \cdot \rangle} \right) = d_U f \otimes s e^{\sqrt{-1} \langle \mu, \cdot \rangle}.$$

Thus, acting on $C^\infty(U) \, \widehat{\otimes} \, \Lambda\left((\mathbb{C}^N)^*\right) \otimes \mathbb{C} e^{\sqrt{-1} \langle \mu, \cdot \rangle}$, we have

$$(B.40) \qquad d^M = \left(I \, \widehat{\otimes} \, \sqrt{-1} \, (\mu \wedge) \otimes I\right) + \left(d_U \, \widehat{\otimes} \, I \otimes I\right).$$

Using the isomorphism

$$(B.41) \qquad C^\infty(U; \Lambda(E^*|_U)) \cong C^\infty(U) \otimes \Lambda\left((\mathbb{C}^N)^*\right),$$

we see that acting on $C^\infty(U; \Lambda(E^*|_U)) \otimes \mathbb{C} e^{\sqrt{-1} \langle \mu, \cdot \rangle}$, one has

$$(B.42) \qquad d^M = \left( \sqrt{-1} \, (\mu \wedge) + d_U \right) \otimes I.$$

However, in terms of our trivializations, $\nabla^{\Lambda(E^*|_U)} = d_U$. The theorem follows. ∎

**Theorem B.14 :** *With respect to the orthogonal decomposition* (B.34), *the superconnection* $(d^M)^*$ *splits as*

$$(B.43) \qquad (d^M)^* = \bigoplus_{\mu \in C^\infty(U;\Lambda^*)} \left( -\sqrt{-1}\, i(\mu) + \nabla^{\Lambda(E^*|_U)} + \Lambda\left(\omega_{E^*|_U}\right) \right) \otimes I.$$

**Proof :** Let $h^{\mathbb{C}}$ be the standard Hermitian metric on $\mathbb{C}$. The restriction of $h^W$ to $L^2(U; \Lambda(E^*|_U)) \otimes \mathbb{C} e^{\sqrt{-1}\langle \mu, \cdot \rangle}$ is $h^{\Lambda(E^*|_U)} \otimes h^{\mathbb{C}}$. Thus it is enough to compute the adjoint of $\sqrt{-1}\,(\mu \wedge) + \nabla^{\Lambda(E^*|_U)}$, acting on $L^2(U; \Lambda(E^*|_U))$, with respect to the Hermitian metric $h^{\Lambda(E^*|_U)}$. Clearly the adjoint of $(\mu \wedge)$ is $i(\mu)$. The adjoint of $\nabla^{\Lambda(E^*|_U)}$ is

$$(B.44)\ \nabla^{\Lambda(E^*|_U)} + \left(h^{\Lambda(E^*|_U)}\right)^{-1} \left(\nabla^{\Lambda(E^*|_U)} h^{\Lambda(E^*|_U)}\right) = \nabla^{\Lambda(E^*|_U)} + \omega_{\Lambda(E^*|_U)}$$

$$= \nabla^{\Lambda(E^*|_U)} + \Lambda\left(\omega_{E^*|_U}\right).$$

The theorem follows. ∎

We can now remove the restriction to the open set $U$, and work globally on $B$. Put

$$(B.45) \qquad \widehat{c}(\mu) = (\mu \wedge) + i(\mu).$$

For $t > 0$ and $\mu \in \Lambda^*$, put

$$(B.46) \qquad D_{t,\mu} = -\sqrt{-1}\,\frac{\sqrt{t}}{2}\,\widehat{c}(\mu) + \frac{1}{2}\,\Lambda\left(\omega_{E^*}\right).$$

Defining $D_t$ as in (A.38), Theorems B.13 and B.14 give that

$$(B.47) \qquad D_t = \bigoplus_{\mu \in \Lambda^*} D_{t,\mu} \otimes I.$$

Then defining $f\left(C'_t, h^W\right)$ as in (A.39) and $f^\wedge\left(C'_t, h^W\right)$ as in (A.40), we have

$$(B.48) \qquad f\left(C'_t, h^W\right) = \sum_{\mu \in \Lambda^*} (2i\pi)^{1/2} \varphi \operatorname{Tr}_s\left[f(D_{t,\mu})\right],$$

$$f^\wedge\left(C'_t, h^W\right) = \sum_{\mu \in \Lambda^*} \varphi \operatorname{Tr}_s\left[\frac{N}{2} f'(D_{t,\mu})\right],$$

where the supertraces on the right-hand-side of (B.48) are finite-dimensional supertraces on $\Lambda(E^*)$. Applying the results of the previous section with $\mathcal{V} = E^*$, we see from (B.48) and Definition B.8 that

$$(B.49) \qquad f\left(C'_t, h^W\right) = \sum_{\mu \in \Lambda^*} \mu^* \widetilde{\delta}_t,$$

$$\frac{1}{t} f^\wedge\left(C'_t, h^W\right) = \sum_{\mu \in \Lambda^*} \mu^* \widetilde{\epsilon}_t.$$

**Theorem B.15 :** *One has*

$$(B.50) \qquad \frac{\partial}{\partial t} f\left(C'_t, h^W\right) = \frac{1}{t} df^\wedge\left(C'_t, h^W\right).$$

**Proof :** This follows from Theorem A.16, or more directly from equation (B.15). ∎

**Theorem B.16 :** *If $K$ is a compact subset of $B$ then there is a constant $c > 0$ such that on $K$, as $t \to \infty$, one has*

$$(B.51) \qquad f\left(C'_t, h^W\right) = f\left(\nabla^{\Lambda(E^*)}, h^{\Lambda(E^*)}\right) + \mathcal{O}(e^{-ct}),$$

$$f^\wedge\left(C'_t, h^W\right) = -\frac{1}{2} + \mathcal{O}(e^{-ct}) \quad \text{if } \dim(E) = 1$$
$$= \mathcal{O}(e^{-ct}) \qquad \text{if } \dim(E) > 1.$$

**Proof :** The contribution of the terms in (B.49) with $\mu \neq 0$ is exponentially small in $t$. Thus it suffices to look at the $\mu = 0$ term. As

$$(B.52) \qquad (2i\pi)^{1/2} \varphi \operatorname{Tr}_s\left[f\left(D_{t,0}\right)\right] = f\left(\nabla^{\Lambda(E^*)}, h^{\Lambda(E^*)}\right),$$

the first line of (B.51) follows. Now

$$(B.53) \qquad \varphi \operatorname{Tr}_s\left[\frac{N}{2} f'\left(D_{t,0}\right)\right] = \varphi \operatorname{Tr}_s\left[\frac{N}{2} f'\left(\frac{1}{2} \Lambda\left(\omega_{E^*}\right)\right)\right].$$

By Lemma 1.3, this is the same as

$$(B.54) \qquad \operatorname{Tr}_s\left[\frac{N}{2} f'(0)\right] = \frac{1}{2} \sum_{p=0}^N (-1)^p\, p \binom{N}{p},$$

from which the rest of Theorem B.16 follows. ∎

**Remark B.17 :** Theorem B.16 is the same as Theorem 2.17, which was proved in terms of Berezin integrals.

Combining Theorem B.11 and (B.49), we see that the form $\sum_{\mu \in \Lambda} \mu^* \delta_t$ of Section 2d is essentially the same as the form $f\left(C'_t, h^W\right)$ and the form $\sum_{\mu \in \Lambda} \mu^* \epsilon_t$ of Section 2d is essentially the same as the form $\frac{1}{t} f^\wedge \left(C'_t, h^W\right)$. Thus the analysis of Section 2d simply consists of explicitly verifying the general properties of the forms $f\left(C'_t, h^W\right)$ and $f^\wedge \left(C'_t, h^W\right)$ in this special case of a torus fibration.